\documentclass[12pt]{amsart}
\usepackage{amssymb}  
\usepackage[all]{xy}
\def\charic{{\operatorname{char}\,}} 
\def\dv{\!\mid\!}

\def\ratto{\,\,-\,\!\to}

\def\Gm{\mathbb G_m}
\def\PP{\mathbb P}
\def\F{\mathbb F}

\def\Div{\operatorname{Div}}
\def\Tr{\operatorname{Tr}}
\def\GL{\operatorname{GL}}
\def\tors{\operatorname{tors}}

\def\Gal{\operatorname{Gal}}
\def\Br{\operatorname{Br}}
\def\PGL{\operatorname{PGL}}
\def\Aut{\operatorname{Aut}}
\def\kbar{\overline{k}}
\def\SL{\operatorname{SL}}
\def\Diag{\operatorname{Diag}}

\def\ord{\operatorname{ord}}
\def\div{\operatorname{div}}
\def\adj{\operatorname{adj}}
\def\summ{\sigma}
\def\isom{{\, \cong \,}}
\newcommand{\Q}{{\mathbb Q}}

\newcommand{\Z}{{\mathbb Z}}

\def\ra{{\longrightarrow}}
\def\Mat{\operatorname{Mat}}
\def\End{\operatorname{End}}
\def\Ob{\operatorname{Ob}}
\def\oh{{\mathcal O}}
\def\L{{\mathcal L}}

\def\la{{\lambda}}
\def\x{{\mathbf x}}
\def\u{{\mathbf u}}
\def\h{{\mathbf h}}
\def\v{{\mathbf v}}
\def\t{{\mathbf t}}
\def\e{{\mathbf e}}
\def\f{{\mathbf f}}

\def\a{{\mathbf a}}
\newcommand{\Sha}{\mbox{\wncyr Sh}}
\newfont{\wncyr}{wncyr10 at 12pt}

\newenvironment{Proof}{\par\noindent{\sc Proof:}}%
                      {\hspace*{\fill}\nobreak$\Box$\par\medskip}
\newenvironment{ProofOf}[1]{\par\noindent{\sc Proof of #1:}}%
                      {\hspace*{\fill}\nobreak$\Box$\par\medskip}

\newtheorem{Proposition}{Proposition}[section]
\newtheorem{Theorem}[Proposition]{Theorem}
\newtheorem{Lemma}[Proposition]{Lemma}
\newtheorem{Corollary}[Proposition]{Corollary}

\theoremstyle{definition}

\newtheorem{Definition}[Proposition]{Definition}
\newtheorem{Remark}[Proposition]{Remark}
\newtheorem{Algorithm}[Proposition]{Algorithm}

\begin{document}

\date{23rd November 2007}
\title[Finding rational points on elliptic curves]
{Finding rational points on elliptic curves 
using 6-descent and 12-descent}

\author{Tom Fisher}
\address{University of Cambridge,
        DPMMS, Centre for Mathematical Sciences,
         Wilberforce Road, Cambridge CB3 0WB, UK}
\email{T.A.Fisher@dpmms.cam.ac.uk}

 \begin{abstract}
We explain how recent work on $3$-descent and $4$-descent for elliptic
curves over $\Q$ can be combined to search for generators of the 
Mordell-Weil group of large height. As an application we show that 
every elliptic curve of prime conductor in the Stein-Watkins 
database has rank at least as large as predicted by the conjecture 
of Birch and Swinnerton-Dyer.
 \end{abstract}

\maketitle

\section{Introduction}

Let $E$ be an elliptic curve over $\Q$. An 
$n$-descent calculation on $E$ provides us with 
$n$-covering curves $\pi_\alpha : C_\alpha \to E$ for 
$\alpha$ running over a finite indexing set $A$, with 
the property that 
$$ \bigcup_{\alpha \in A} \pi_\alpha ( C_\alpha(\Q) ) = E(\Q). $$
The usual choice of indexing set $A$ is the $n$-Selmer
group $S^{(n)}(E/\Q)$ which sits in a short exact sequence
$$ 0 \to E(\Q) / nE(\Q) \stackrel{\delta}{\to} 
S^{(n)}(E/\Q) \to \Sha(E/\Q)[n] \to 0. $$
Given $\alpha \in S^{(n)}(E/\Q)$ there are two possibilities:
either $\pi_{\alpha}(C_\alpha(\Q))$ is a coset of 
$nE(\Q)$ in $E(\Q)$, in which case $\alpha$ is the image of this
coset by $\delta$, or $C_\alpha(\Q)$ is empty, 
in which case $\alpha$ maps to a non-trivial element 
of the Tate-Shafarevich group $\Sha(E/\Q)$.

It has long been known that $n$-descent can help in the
search for generators of the Mordell-Weil group $E(\Q)$. Indeed the
theory of heights (see for example \cite{Hindry-Silverman}) 
suggests that if we
write our $n$-coverings as curves of degree $n$ with small
coefficients, then a point of (logarithmic) height $h$ on $E(\Q)$
should come from a point of height approximately $h/(2n)$ 
on $C_\alpha(\Q)$ for suitable $\alpha$. This is not a precise statement
(the height is only bounded up to the addition of a constant whose 
behaviour with respect to $n$ is unknown)
but the idea seems to work well in practice.

We would therefore like to perform $n$-descent calculations
for $n$ as large as possible. Until recently $n$-descent has 
only been practical for general\footnote{{\em i.e.} we make no 
assumption on the Galois module structure of $E[n]$.}  
elliptic curves in the case $n=2$. 
Methods for $4$-descent and $8$-descent have been developed in
the PhD theses of Siksek~\cite{Siksek}, Womack~\cite{Womack} and 
Stamminger~\cite{Stamminger}. 
Joint work of the author with 
Cremona, O'Neil, Simon and Stoll \cite{descsum} has now made 
$3$-descent practical, and in a few preliminary examples also $5$-descent. 
The algorithms for $2$-descent, $3$-descent and $4$-descent have been 
contributed to the computer algebra system Magma \cite{Magma}.

The natural question therefore arises as to how we can combine
an $m$-covering and $n$-covering for $m$ and $n$ coprime to 
give an $mn$-covering. At the level of abelian groups it is trivial
that 
$$ S^{(mn)}(E/\Q) \isom S^{(m)}(E/\Q) \times S^{(n)}(E/\Q). $$
However, if we are to represent the Selmer group elements
as covering curves, then it is not so clear how one should
proceed.

Suppose we are given an $m$-covering $\pi_m : C_m \to E$ and an 
$n$-covering $\pi_n : C_n \to E$. Then the curves $C_m$ and $C_n$ 
are torsors under $E$, and an $mn$-covering is given by 
$$ C_{mn} = \frac{C_m \times C_n}{E} $$
where we quotient out by the diagonal action of $E$. 
An alternative would be to take fibre product 
$$ C_{mn} = C_m \times_{E} C_n $$
with respect to the covering maps $\pi_m$ and $\pi_n$.
As far as we can see, neither of these constructions is suitable
for practical computation. 
We have therefore taken a different approach based on 
representations of the Heisenberg group. 

Unfortunately our approach does not work for arbitrary coprime integers
$m$ and $n$, but only when each of $m$ and $n$ is plus or minus 
a square modulo the other. 
This includes the case of consecutive integers. 
In this case, we specify
an embedding of $E$ in $\PP(\Mat_{n,n+1})$ as a curve of degree
$n(n+1)$, in such a way that when $E$ acts on itself by translation,
the $n$-torsion points act as left multiplication by $n \times n$ matrices, 
and the $(n+1)$-torsion points act as right 
multiplication by $(n+1) \times (n+1)$ matrices.
We can then twist $E$ by a pair of
cocycles taking values in $E[n]$ and $E[n+1]$ to obtain the
required $n(n+1)$-covering $C_{n(n+1)}$ as a curve in $\PP(\Mat_{n,n+1})$.
Moreover, it turns out that the covering map $C_{n(n+1)} \to C_{n+1}$ 
is defined by the $n \times n$ minors.

We give a precise statement of these results in \S\ref{sec:comptwists}.
We employ two different methods of
proof. The first, described in \S\ref{sec:heis}, uses representations of 
the Heisenberg group, and leads to results for arbitrary $n$. The second, 
described in \S\S\ref{invthy},\ref{sec:covmat}, 
uses the invariant theory of binary quartics and ternary cubics,
and gives practical formulae specific to $6$-descent and $12$-descent. 

In \S\ref{implement} we give some details of our implementation of 
6-descent and 12-descent in Magma.
Using 12-descent, we now expect to be able to find
rational points on an elliptic curve over $\Q$ 
up to logarithmic height 600 (provided the 
coefficients of the original elliptic curve are not too large).
The main bottleneck comes in the $3$-descent, where we must
compute the class group and units of each number field
generated by the co-ordinates of a 3-torsion point of $E$.
(There is usually just one such field, and it has degree $8$.)
Fortunately, since our final answer comes in the form of a
rational point, there is no need to perform these intermediate
calculations rigorously.

Stein and Watkins \cite{SteinWatkins} have constructed a database 
of elliptic curves that is expected to contain most elliptic 
curves over $\Q$ of prime conductor $N$ with $N < 10^{10}$. 
We are able to show that every curve in their 
database (of prime conductor) has rank at least as large as 
predicted by the conjecture of Birch and Swinnerton-Dyer.
Prior to our involvement, this had been reduced by 
Cremona and Watkins to a list of 35 curves of analytic rank 2 
for which one generator of small height (less than $34$) was known, but 
a second generator of large height (greater than $220$) remained 
to be found. In each case we were able to find the second generator
using either 6-descent or 12-descent. 

We give two numerical examples in \S\ref{sec:num}. In the first we use 
6-descent to find a pair of non-zero integers $x$ and $y$ for which both 
$$x^2 - 809xy+y^2 \quad \text{ and } \quad x^2 + 809xy+y^2$$ are
squares. We find that the smallest solution is given by 
$x$ and $y$ with $534$ and $537$ decimal digits respectively.
Our second example is the last in the list of 35 curves mentioned
above. In this case we use 12-descent to find a generator of 
height $642.63.$

\section*{Acknowledgments}

I would like to thank Steve Donnelly for sharing his
initial thoughts on this problem, and Mark Watkins for
suggesting suitable test data, including the examples 
in \S\ref{sec:num}. All computer calculations in support
of this work were performed using Magma \cite{Magma}.

\section{Computing twists}
\label{sec:comptwists}

Let $k$ be a field of characteristic zero, 
with algebraic closure $\kbar$.
We fix an elliptic curve $E$ over $k$ with identity $\oh$. 
Let $n \ge 2$ be an integer. 
A base diagram $[E \to \PP^{n-1}]$ of level $n$
is a morphism defined over $k$ determined by the complete linear 
system $|n.\oh|$. Thus any two base diagrams differ by an 
element of $\PGL_n(k)$.

More generally we consider diagrams $[C \to S]$ where $C$ is
a torsor under $E$ and $S$ is a variety  (both defined over $k$) 
and the map $C \to S$ is a morphism defined over $k$. Two such diagrams
$[\phi_1 : C_1 \to S_1]$ and $[\phi_2: C_2 \to S_2]$ are isomorphic
if there is an isomorphism of torsors $\alpha : C_1 \isom C_2$ 
and an isomorphism of varieties $\beta : S_1 \isom S_2$
satisfying $\phi_2 \circ \alpha = \beta \circ \phi_1$. 
We define a Brauer-Severi diagram $[C \to S]$ to be a twist of the
base diagram. Then $S$ is a Brauer-Severi variety, and the 
morphism $C \to S$ is that determined by a complete 
linear system $|D|$, where the divisor $D$ is  
linearly equivalent to all its Galois conjugates, 
but need not itself be defined over $k$.

We recall from \cite[Paper I, \S1.3]{descsum} that the Brauer-Severi diagrams 
are parametrised, as twists of a fixed base diagram, by the
Galois cohomology group $H^1(k,E[n])$. Moreover there is an 
obstruction map 
$$ \Ob_n : H^1(k,E[n]) \to \Br(k)[n] $$
taking the class of $[C \to S]$ to the class of $[S]$. 
In general this map is not a group homomorphism. We are interested
in the elements of $H^1(k,E[n])$ with trivial obstruction,
equivalently those that are represented by diagrams of the form 
$[\phi_n : C \to \PP^{n-1}]$. With the convention that
points of $\PP^{n-1}$ are written as column vectors, 
we define the ``character'' associated to $\phi_n$ to be the 
unique morphism of $k$-group schemes 
$\chi_n : E[n] \to \PGL_n$ such that
$$ \phi_n (T + P) = \chi_n(T) \phi_n (P) $$
for all $T \in E[n](\kbar)$ and $P \in C(\kbar)$.

In \S\ref{sec:heis} we use representations of the Heisenberg group 
to prove

\begin{Theorem}
\label{thm1}
Let $m$ and $n$ be coprime integers satisfying
\begin{equation}
\label{mnab} 
u^2 n \equiv \pm 1 \pmod{m} \quad \text{ and } \quad
v^2 m \equiv \pm 1 \pmod{n}, 
\end{equation}
for some integers $u$ and $v$. 
Let $\chi_m$, $\chi_n$ and $\chi_{mn}$ be the characters associated
to base diagrams of level $m$, $n$ and $mn$.
Then there is a morphism of $k$-group schemes
$$ \Xi: \PGL_m \times \PGL_n \to \PGL_{mn} $$
such that
 $$\Xi( \chi_m(S), \chi_n(T)) = \chi_{mn}(u S + v T)$$ 
for all $S \in E[m](\kbar)$ and $T \in E[n](\kbar)$. 
\end{Theorem}

Let $[C \to \PP^{n-1}]$ be a twist of the base diagram 
$[E \to \PP^{n-1}]$.
By definition this means that there is a commutative diagram
$$ \xymatrix{ C \ar[r] \ar[d]_\alpha & \PP^{n-1} \ar[d]^\beta \\ 
E \ar[r] & \PP^{n-1} } $$
where $\alpha$ and $\beta$ are isomorphisms defined over $\kbar$.
We say that the matrix $B_n \in \PGL_n(\kbar)$ representing $\beta$ is 
a flex matrix for $[C \to \PP^{n-1}]$.
In \S\ref{compfm} we give some algorithms for computing flex matrices,
starting from equations for $C$. Conversely, it is clear that we can recover 
equations for $C$ from a flex matrix (by starting with equations for
$E$ and making the relevant substitution).

Our next result explains how Theorem~\ref{thm1} is used to construct 
an $mn$-covering from an $m$-covering and an $n$-covering.

\begin{Proposition}
In the setting of Theorem~\ref{thm1}, suppose that 
$\xi \in H^1(k,E[m])$ and $\eta \in H^1(k,E[n])$ 
are represented by Brauer-Severi diagrams $[C_m \to \PP^{m-1}]$ and 
$[C_n \to \PP^{n-1}]$ with flex matrices $B_m$ and $B_n$.
Then $u \xi + v \eta \in H^1(k,E[mn])$ is represented by a 
Brauer-Severi diagram $[C_{mn} \to \PP^{mn-1}]$ with flex matrix 
$B_{mn} = \Xi(B_m,B_n)$. 
\end{Proposition}

\begin{Proof}
Let $\sigma \in \Gal(\kbar/k)$. 
Since $\sigma(B_m) B_m^{-1} = \chi_m(\xi_\sigma)$, 
$\sigma(B_n) B_n^{-1} = \chi_n(\eta_\sigma)$ and $\Xi$ is
defined over $k$, we have
$$ \begin{array}{rcl}
\sigma(B_{mn}) B_{mn}^{-1} & = & \Xi(\sigma(B_m) B_m^{-1}, 
\sigma(B_n) B_n^{-1}) \\
& = & \Xi( \chi_m(\xi_\sigma), \chi_n(\eta_\sigma)) \\
& = & \chi_{mn} ( u  \xi_\sigma + v \eta_\sigma).
\end{array} $$
Thus $B_{mn}$ represents a change of co-ordinates on $\PP^{mn-1}$
taking $E$ to a curve $C_{mn}$ defined over $k$. 
Then $[C_{mn} \to \PP^{mn-1}]$ is the twist of $[E \to \PP^{mn-1}]$
by $u \xi + v \eta \in H^1(k,E[mn])$.
\end{Proof}

\begin{Remark} The proposition
shows that if $\xi$ and $\eta$ have trivial obstruction, then 
so does $u \xi + v \eta$. In fact, standard properties of the
obstruction map (see \cite[Paper I]{descsum}, \cite{Cathy}, \cite{Zarhin}) 
already show that
\begin{equation}
\label{obcalc}
 \Ob_{mn}( u \xi+ v \eta) 
  = u^2 n \Ob_m(\xi) + v^2 m \Ob_n(\eta) 
\end{equation}
for $m$ and $n$ coprime. Interestingly, 
the hypothesis~(\ref{mnab}) of 
Theorem~\ref{thm1} is that the coefficients on the 
right hand side of~(\ref{obcalc}) are $\pm 1$.
\end{Remark}

We are mainly interested in the case $m$ and $n$ are consecutive 
integers.
We can therefore either take $u=v=1$ in Theorem~\ref{thm1}, 
or use the refined version of the theorem 
we give next.

We write $\PP(\Mat_{a,b})$ for the projective space of dimension 
$ab-1$ formed from the vector space of $a \times b$ matrices. 
Taking $n \times n$ minors defines a rational map 
$$ \mu : \PP (\Mat_{n,n+1}) \ratto \PP^{n}\, ; \quad 
A \mapsto ( \ldots : (-1)^i \det(A^{\{i\}}) : \ldots ) $$
where $A^{\{i\}}$ is $A$ with the $i$th column deleted.
The following refinement of Theorem~\ref{thm1} is proved
alongside the original theorem in \S\ref{sec:heis}.

\begin{Theorem}
\label{thm2}
Let $[\phi_n : E \to \PP^{n-1}]$ and $[\phi_{n+1}: E \to \PP^n]$ 
be base diagrams of levels $n$ and $n+1$, with associated characters
$\chi_n$ and $\chi_{n+1}$. Then there is a base diagram 
$$ [\phi_{n,n+1} : E \to \PP(\Mat_{n,n+1})] $$
of level $n(n+1)$, with associated character $\chi_{n,n+1}$ given by 
$$ \chi_{n,n+1}(S+T) : A \mapsto \chi_{n}(S) A \chi_{n+1}(T)$$
for all $S \in E[n](\kbar)$ and $T \in E[n+1](\kbar)$.
Moreover if $[n] : E \to E$ is the multiplication-by-$n$ map 
then the diagram
$$ \xymatrix{ E \ar[rr]^-{\phi_{n,n+1}} \ar[d]_{[n]} & & 
 \PP (\Mat_{n,n+1}) \ar@{-->}[d]^\mu \\
E \ar[rr]^-{\phi_{n+1}} & & \PP^{n}. } $$
commutes.
\end{Theorem}

We obtain equations for an $n(n+1)$-covering from an $n$-covering and 
an $(n+1)$-covering, by first finding the
base diagram $\phi_{n,n+1}$ of Theorem~\ref{thm2}, and then twisting
by the flex matrices $B_n$ and $B_{n+1}$. These twists may
be performed one after the other. In fact twisting by
$B_n$ first gives the following generalisation of Theorem~\ref{thm2}.

\begin{Theorem} 
\label{thm3}
Let $[\phi_n : C \stackrel{|D|}{\ra} \PP^{n-1}]$ be a Brauer-Severi diagram
and $[\phi_{n+1} : E \to \PP^n]$ a base diagram, of levels $n$ 
and $n+1$, with associated characters $\chi_n$  and $\chi_{n+1}$.
Then there is a Brauer-Severi diagram
$$ [\phi_{n,n+1}: C \stackrel{|(n+1)D|}{\ra} \PP(\Mat_{n,n+1})] $$
of level $n(n+1)$, with 
associated character $\chi_{n,n+1}$ given by 
$$ \chi_{n,n+1}(S+T) : A \mapsto \chi_{n}(S) A \chi_{n+1}(T)$$
for all $S \in E[n](\kbar)$ and $T \in E[n+1](\kbar)$.
Moreover if $\pi : C \to E$ is the $n$-covering map then the diagram
$$ \xymatrix{ C \ar[rr]^-{\phi_{n,n+1}} \ar[d]_{\pi} & & 
 \PP (\Mat_{n,n+1}) \ar@{-->}[d]^\mu \\
E \ar[rr]^-{\phi_{n+1}} & & \PP^{n} } $$
commutes.
\end{Theorem}
\begin{Proof} 
Theorem~\ref{thm2} is the special case where $(C,[D]) = (E,[n.\oh])$. 
We twist by $B_n$ to obtain the general result.
\end{Proof}

The advantage of using Theorem~\ref{thm3}  (instead of Theorem~\ref{thm2})
is that we then only need to twist by $B_{n+1}$ (rather than both $B_n$
and $B_{n+1}$) to obtain the desired $n(n+1)$-covering.
In \S\ref{sec:covmat} we use invariant theory to give an alternative
proof of Theorem~\ref{thm3} in the cases $n=2,3$.
In particular we obtain explicit formulae for $\phi_{2,3}$ and
$\phi_{3,4}$ which we then use in 
our implementations of 6-descent and 12-descent.

\section{The Heisenberg group}
\label{sec:heis}

We continue to work over a field $k$ of characteristic $0$.
Let $E$ be an elliptic curve with identity $\oh$, and $n \ge 2$ an integer. 
We recall that if $D$ is a divisor on $E$
of degree $n$ then the Riemann-Roch space $\L(D)$ has dimension $n$.  
Let $V_n = \L(n.\oh)^*$. Then there is a ``co-ordinate free'' 
base diagram $[E \to \PP(V_n)]$ with associated character
$\chi_n : E[n] \to \PGL(V_n)$.

\begin{Definition}
(i) The theta group $\Theta_n$ is the inverse image of $\chi_n(E[n])$
in $\GL(V_n)$. It sits in a commutative diagram of $k$-group
schemes with exact rows:
$$ \xymatrix{ 0 \ar[r] & \Gm \ar[r] \ar@{=}[d] 
& \Theta_n \ar[r] \ar[d] & E[n] \ar[r] \ar[d]^{\chi_n} & 0 \\
 0 \ar[r] & \Gm \ar[r] & \GL(V_n) \ar[r] & \PGL(V_n) \ar[r] & 0 \rlap{.}} $$
(ii) The Heisenberg group $H_n$ is the inverse image of $\chi_n(E[n])$
in $\SL(V_n)$. It sits in a commutative diagram of $k$-group
schemes with exact rows:
$$ \xymatrix{ 0 \ar[r] & \mu_n \ar[r] \ar@{=}[d] 
& H_n \ar[r] \ar[d] & E[n] \ar[r] \ar[d]^{\chi_n} & 0 \\
 0 \ar[r] & \mu_n \ar[r] & \SL(V_n) \ar[r] & \PGL(V_n) \ar[r] & 0 \rlap{.}} $$
\end{Definition}

\begin{Remark}
\label{RemH}
It is well known (see e.g. \cite{Hulek}) 
that over $k=\kbar$ we may choose a basis for
$V_n$ such that $\Theta_n$ is generated by 
$$ \begin{pmatrix}
    1   &    0   &    0    & \cdots &      0     \\
    0   & \zeta_n  &    0    & \cdots &      0     \\
    0   &    0   & \zeta_n^2 & \cdots &      0     \\ 
 \vdots & \vdots & \vdots  &        &   \vdots   \\
    0   &   0    &    0    & \cdots & \zeta_n^{n-1}  
\end{pmatrix}, \quad  
    \begin{pmatrix}
   0   &    0   & \cdots &   0    &      1     \\ 
    1   &    0   & \cdots &   0    &      0     \\
    0   &    1   & \cdots &   0    &      0     \\ 
 \vdots & \vdots &        & \vdots &   \vdots   \\
    0   &   0    & \cdots &   1    &      0       
\end{pmatrix}, $$
and the scalar matrices. (Here $\zeta_n \in k$ is a primitive
$n$th root of unity.) In particular, taking commutators in 
$\Theta_n$ defines a non-degenerate pairing 
$e_n : E[n] \times E[n] \to \mu_n$,
which turns out to be the Weil pairing.
\end{Remark}

The group $H_n(\kbar)$ is a non-abelian group of order $n^3$
with centre $\mu_n(\kbar)$.

\begin{Definition}
A representation $\rho : H_n(\kbar) \to \GL_d(\kbar)$ has central
character $[r]$ if $\rho(\zeta) = \zeta^r I_d$ 
for all $\zeta \in \mu_n(\kbar)$.
\end{Definition}

\begin{Lemma} 
\label{LemmaR}
Let $r$ be an integer coprime to $n$. Then \\
(i) Every representation of $H_n(\kbar)$ with central character $[r]$
is a direct sum of irreducible $n$-dimensional representations. \\
(ii) Up to equivalence, there is a unique $n$-dimensional representation
of $H_n(\kbar)$ with central character $[r]$. 
\end{Lemma}
\begin{Proof}
Let $\rho$ be a representation of $H_n(\kbar)$ with central 
character $[r]$. If $\sigma, \tau \in H_n(\kbar)$ lift 
$S,T \in E[n](\kbar)$ then 
$ \rho(\sigma) \rho(\tau) \rho(\sigma)^{-1}  = e_n(S,T)^r \rho(\tau) $
and so $\Tr(\rho(\tau)) = e_n(S,T)^r \Tr(\rho(\tau))$. 
Since $r$ is coprime to $n$, the character of $\rho$ 
vanishes outside the centre of $H_n(\kbar)$.
Then orthogonality relations in the character table show that $\rho$ is
irreducible if and only if it has dimension $n$. 
This proves (i) and the uniqueness in (ii). 
Existence is clear in the case $r=1$. 
In general we take the $r$th tensor power and apply (i). 
\end{Proof}

To prove Theorem~\ref{thm1} we need to construct a morphism
$\Xi$ that is defined over $k$. We therefore study
representations $\rho : H_n(\kbar) \to \GL_d(\kbar)$ that are Galois
equivariant, equivalently those that induce a morphism of 
$k$-group schemes $H_n \to \GL_d$. We call these representations of $H_n$.

By construction, $V_n$ is an $n$-dimensional representation of $H_n$
with central character $[1]$. So by Lemma~\ref{LemmaR}(i) it is
irreducible. We might hope to construct other irreducible $n$-dimensional
representations of $H_n$ by any one of the following standard methods.

\renewcommand{\labelenumi}{(\roman{enumi})}
\begin{enumerate}
\item Take a subspace or quotient of a tensor power of $V_n$. 
\item Replace $V_n$ by one of its Galois conjugates.
\item Precompose $\rho : H_n \to \GL(V_n)$ with an automorphism of $H_n$.
\end{enumerate}

We see no way of using (i) in the proof of Theorem~\ref{thm1}, other
than in the case of $V_n^* = \wedge^{n-1} V_n$. (The problem is 
that there is no analogue of Lemma~\ref{LemmaR}(i) with $H_n$ 
replaced by $\GL_n$.) Our restriction to Galois equivariant 
representations rules out the use of (ii). 

To use (iii) we must first
describe the (Galois equivariant) automorphisms of $H_n$. 
Each automorphism of $H_n$ induces an automorphism of $E[n]$. We 
may identify $\Aut(E[n])$ as the centraliser of the image of Galois
in $\Aut_{\kbar}(E[n]) \isom \GL_2(\Z/n\Z)$. So in general the
only automorphisms of $E[n]$ are the maps $T \mapsto aT$ for 
$a \in (\Z/n\Z)^\times$. We show that each of these maps lifts to 
an automorphism of $H_n$. (Without our insistence on Galois 
equivariance, this would be trivial.)

\begin{Lemma}
\label{LemmaLambda}
For each $a \in (\Z/n\Z)^\times$ there is a morphism of $k$-group
schemes $\psi_a : H_n \to H_n$ making the diagram 
$$ \xymatrix{ 0 \ar[r] & \mu_n \ar[r] \ar[d]^{a^2} 
& H_n \ar[r] \ar[d]^{\psi_a} & E[n] \ar[r] \ar[d]^a & 0 \\
 0 \ar[r] & \mu_n \ar[r] & H_n \ar[r] & E[n] \ar[r] & 0 } $$
commute.
\end{Lemma}

\begin{Proof} The $[-1]$-map on $E$ lifts to a matrix $\iota \in \GL_n(k)$. 
Conjugation by $\iota$ gives the map $\psi_{-1}$. In general
we first define
$$ \lambda : H_n \to \mu_n \, ; \quad h \mapsto \iota h \iota^{-1} h $$
and then put
$$ \psi_a : H_n \to H_n \, ; \quad h \mapsto \lambda(h)^{a(a-1)/2} h^a. $$
Since $\lambda(xy) = xy (\iota x \iota^{-1})(\iota y \iota^{-1}) =
 \lambda(x) \lambda(y) x y x^{-1} y^{-1}$ 
for all $x,y \in H_n$, it is easy to check that 
$\psi_a$ is a group homomorphism. Galois equivariance
is clear from the construction.
\end{Proof}

Let $V_n^{(a)}$ be the representation of $H_n$ given by
$$ H_n \times V_n \to V_n \, ; \quad (h,v) \mapsto \psi_a(h)v. $$
It is an irreducible $n$-dimensional 
representation of $H_n$ with central character $[a^2]$. 
Taking these representations and their duals we obtain
all the representations of Lemma~\ref{LemmaR}(ii) with
$r \equiv \pm a^2 \pmod{n}$. 

We write $\tau_P : E \to E$ for translation by $P \in E(\kbar)$. 
We recall from \cite[\S23]{Mumford}, (see also \cite[Paper I]{descsum}), 
that the theta group $\Theta_n$ may be described as pairs 
$$ \Theta_n(\kbar) = \{ \, (f,T) \in  \kbar(E)^\times \times E[n](\kbar) 
\mid \div(f) = \tau_T^*(n. \oh)- n.\oh \, \} $$
with group law
\begin{equation}
\label{gplaw}
 (f,S) * (g,T) = ((\tau_T^*f)g,S+T). 
\end{equation} 
There is a natural action of 
$\Theta_n$ on $V_{nN}^* = \L(nN.\oh)$  given by 
\begin{equation}
\label{HactV*}
(f,T) : \, h \mapsto \tau_{-T}^* ( h/f^N ). 
\end{equation}

We use this notation to relate the Heisenberg groups $H_n$ for  
different levels $n$.

\begin{Proposition}
\label{PropP}
Let $m$ and $n$ be coprime integers. Then there is an isomorphism
of $k$-group schemes
$$ \begin{array}{rcl} 
H_m \times H_n & \cong & H_{mn} \\ ((f,S),(g,T)) & \mapsto &
 (f^n,S) *(g^m,T). \end{array} $$
\end{Proposition}

For the proof we need two lemmas.

\begin{Lemma}
\label{Lemma1}
Let $M_T \in \Theta_n(\kbar)$ be a lift of $T \in E[n](\kbar)$.
If $T$ has exact order $r$ then 
$$ \det(M_T) = (-1)^{n(n-1)/r} M_T^n. $$
\end{Lemma}
\begin{Proof}
By Remark~\ref{RemH} we know that $M_T$ is similar to 
$$\lambda \Diag(1, \zeta_r, \ldots, \zeta_r^{n-1})$$ 
for some $\lambda \in \kbar^\times$ and 
$\zeta_r$ a primitive $r$th root of unity. Then
$\det(M_T) = \zeta_{r}^{n(n-1)/2} M_T^n = (-1)^{n(n-1)/r} M_T^n$ 
as required.
\end{Proof}

\begin{Lemma}
\label{Lemma2}
Let $N$ be a positive integer. Then there is a commutative
diagram of $k$-group schemes
$$ \xymatrix{ \Theta_n \ar[d]_{\det} \ar[r]^{\alpha} 
& \Theta_{nN} \ar[d]^{\det} \\ \Gm \ar[r]^{N^2} & \Gm \rlap{.} } $$
where $\alpha : (f,T) \mapsto (f^N,T)$.
\end{Lemma}
\begin{Proof} 
It is clear from~(\ref{gplaw}) that $\alpha$ is
a group homomorphism. 
Now let $x \in \Theta_n$ and assume $n$ is odd.
Using Lemma~\ref{Lemma1} we compute
$$ \det (\alpha (x)) = \alpha (x)^{nN} 
  = \alpha (\det x )^{N} = (\det x )^{N^2}. $$
The calculation for $n$ even is similar.
\end{Proof}

\begin{ProofOf}{Proposition~\ref{PropP}} 
Let $(f,S)$ and $(g,T)$ be elements of $\Theta_m(\kbar)$
and $\Theta_n(\kbar)$. The commutator of $(f^n,S)$ 
and $(g^m,T)$ in $\Theta_{mn}(\kbar)$ is both an $m$th root of unity
and an $n$th root of unity, and therefore trivial. 
So there is a morphism of $k$-group schemes
$$ \begin{array}{rcl} 
\Theta_m \times \Theta_n & \to & \Theta_{mn} \\ ((f,S),(g,T)) & \mapsto &
 (f^n,S) *(g^m,T). \end{array} $$
By Lemma~\ref{Lemma2} we can restrict to a map
$H_m \times H_n \to H_{mn}$. Since $m$ and $n$ are coprime, this
restriction is clearly an isomorphism.
\end{ProofOf}

We recall that $H_n$ acts on $\L(nN.\oh)$ as specified in~(\ref{HactV*}).

\begin{Lemma}
\label{charlem}
The $H_n$-invariant subspace of $\L(nN.\oh)$ is trivial unless
$n$ divides $N$, in which case it has dimension $N/n$.
\end{Lemma}
\begin{Proof}
By Lemma~\ref{Lemma2} there is a group homomorphism 
$H_{n}(\kbar) \to H_{nN}(\kbar)$ 
given by $(f,T) \mapsto (f^N,T)$.
By the proof of Lemma~\ref{LemmaR} the character 
of $\L(nN.\oh)$ is trivial outside
the centre of $H_{nN}(\kbar)$. The same is therefore true 
when $\L(nN.\oh)$ is viewed as a representation of $H_{n}(\kbar)$.
We are done by the orthogonality relations in the character table.
\end{Proof}

\begin{Proposition}
\label{PropI}
Let $m$ and $n$ be coprime integers. Suppose that $n \equiv a^2 \pmod{m}$
and $m \equiv b^2 \pmod{n}$ for some integers $a$ and $b$.  
Then there is a $k$-isomorphism of $H_{mn}$-representations 
$$ \pi : V_m^{(a)} \otimes V_n^{(b)} \cong V_{mn}. $$
\end{Proposition}
\begin{Proof}
We recall that $V_m^{(a)}$ is an irreducible $H_m$-representation 
with central character $[a^2]$. Likewise $V_n^{(b)}$ is an 
irreducible $H_n$-representation with central character $[b^2]$.
Then Proposition~\ref{PropP} makes $V_m^{(a)} \otimes V_n^{(b)}$ 
an irreducible $H_{mn}$-representation with central character $[1]$.
Indeed $(\zeta_m,\zeta_n) \in H_m \times H_n$ acts on
$V_m^{(a)} \otimes V_n^{(b)}$ as $\zeta_m^{a^2} \zeta_n^{b^2}$
and on $V_{mn}$ as $\zeta_m^n \zeta_n^m$. The required isomorphism
$\pi$ exists by Lemma~\ref{LemmaR}(ii). Finally, since we work with
Galois equivariant representations, we can
choose an isomorphism $\pi$ that is defined over $k$. 
\end{Proof}

\begin{ProofOf}{Theorem~\ref{thm1}} 
We first treat the case $u^2 n \equiv 1 \pmod{m}$ and 
$v^2 m \equiv 1 \pmod{n}$. Let $a$ and $b$ be inverses
for $u$ and $v$ modulo $m$ and $n$ respectively.
Then the map $\pi : V_m \otimes V_n \to V_{mn}$ constructed in
Proposition~\ref{PropI} satisfies
$$ \pi ( \psi_a(f,S) v_m \otimes \psi_b(g,T) v_n) = 
(( f^n,S) * (g^m,T) ) \pi (v_m \otimes v_n) $$
for all $(f,S) \in H_m$, $(g,T) \in H_n$, $v_m \in V_m$, $v_n \in V_n$.
Passing to $\PP(V_{mn})$ we obtain
$$ \pi ( \chi_m(aS) v_m \otimes \chi_n(bT) v_n ) = \chi_{mn}(S+T) 
\pi (v_m \otimes v_n). $$
Hence 
$$ \pi\circ(\chi_m (aS) \boxtimes \chi_n(bT))=\chi_{mn}(S+T)\circ\pi $$
where $$ \boxtimes : \PGL(V_m) \times \PGL(V_n) \to \PGL(V_m \otimes V_n)$$
is the natural map. The theorem now holds on defining
$$ \begin{array}{rcl}
\Xi : \PGL(V_m) \times \PGL(V_n) & \to & \PGL(V_{mn}) \\
(\alpha,\beta) & \mapsto & \pi \circ (\alpha \boxtimes \beta) \circ \pi^{-1}.
\end{array} $$

In general, if $u^2 n \equiv -1 \pmod{m}$ or 
$v^2 m \equiv -1 \pmod{n}$ then we replace $V_m^{(a)}$ 
or $V_n^{(b)}$ by its dual in Proposition~\ref{PropI},
and the proof carries through as before.
\end{ProofOf}

Next we prove our refined version of the theorem in the case
$m$ and $n$ are consecutive integers.

\medskip

\begin{ProofOf}{Theorem~\ref{thm2}} 
The analogue of Proposition~\ref{PropI} gives a 
$k$-isomorphism of $H_{n(n+1)}$-modules
$$ \pi : V_n \otimes V_{n+1}^* \cong V_{n(n+1)}. $$
Hence there is a base diagram 
\begin{equation}
\label{newbd}
\phi_{n,n+1} : E \to \PP(V_n \otimes V_{n+1}^*) 
\end{equation}
with associated character 
$$ \chi_{n,n+1}(S+T) = \chi_n(S) \boxtimes \chi_{n+1}(T)^*. $$
Picking bases for $V_n$ and $V_{n+1}$, we identify 
$\PP(V_n) = \PP^{n-1}$, $\PP(V_{n+1}) = \PP^{n}$ and $\PP(V_n \otimes
V_{n+1}^*) = \PP(\Mat_{n,n+1})$. Then $\chi_{n,n+1}$ is given by
$$ \chi_{n,n+1}(S+T) : A \mapsto \chi_n(S) A \chi_{n+1} ( \pm T) $$
where the sign $\pm$ is immaterial by the case $a=-1$ of 
Lemma~\ref{LemmaLambda}. This proves the first statement of 
Theorem~\ref{thm2}. 

The base diagram~(\ref{newbd}) is given by a matrix
${\mathfrak A} \in \Mat_{n,n+1}(k(E))$ whose entries are a basis for 
$\L( n(n+1).\oh)$. Let $r$ be the rank of this matrix. 
Then the $r \times r$ minors define a morphism 
$$ \Phi : E \to \PP(\wedge^r V_{n+1}^*) $$
with the property that $\Phi \circ \tau_S = \Phi$ for all
$ S \in E[n](\kbar)$. Hence $\Phi$ factors through $[n] : E \to E$. 
Therefore $n^2 \dv \deg(\Phi^* H)$ where $H$ is the hyperplane section on
$\PP(\wedge^r V_{n+1}^*)$. Since $\deg(\Phi^*H) = r n(n+1)$ 
it follows that $r=n$. Thus there is a commutative diagram
$$ \xymatrix{ E \ar[rr]^-{\phi_{n,n+1}} \ar[d]_{[n]} & & 
 \PP (V_n \otimes V_{n+1}^*) \ar@{-->}[d]^\mu \\
E \ar[rr]^-{\gamma} & & \PP(V_{n+1}) } $$
where the pull back of the hyperplane section by $\gamma$ has degree
$n+1$. We must show that $\gamma = \phi_{n+1}$. 

It is easy to see that $\gamma$ shares with $\phi_{n+1}$ the property
$$ \gamma ( T + P ) = \chi_{n+1} (T) \gamma(P) $$
for all $T \in E[n+1](\kbar)$ and $P \in E(\kbar)$.
Since $V_{n+1}$ is an irreducible representation of $H_{n+1}$
the image of $\gamma$ spans $\PP(V_{n+1})$. So
$\gamma$ is an embedding by a complete linear system.

By Lemma~\ref{charlem} the subspace of $\L (n^2(n+1).\oh)$ fixed by 
$H_n$ has dimension $n+1$. Since the $n \times n$ minors 
of ${\mathfrak A}$ are linearly independent they are 
a basis for this space.
We show in Proposition~\ref{propHfixed} that if $g \in \kbar(E)^\times$ with 
$$\div(g) = (n+1)[n]^* \oh - n^2(n+1).\oh$$ then $g$ is fixed by $H_n$.
Hence we may assume that the first $n \times n$ minor of ${\mathfrak A}$,
viewed as a homogeneous form of degree $n$ in the co-ordinate ring 
of $\PP(V_n \otimes V_{n+1}^*)$, meets $E$ with divisor $(n+1)[n]^* \oh$. 
Then $(n+1).\oh$ is the pull back of a hyperplane section by $\gamma$,
and hence $\gamma$ is a base diagram of level $n+1$.

Since $\gamma$ and $\phi_{n+1}$ are base diagrams of
level $n+1$ they can only differ by an element of $\PGL(V_{n+1})$. 
But they also have the same character $\chi_{n+1}$.
Since the image of $\chi_{n+1}$ is its own centraliser in $\PGL(V_{n+1})$,
it follows that $\gamma = \phi_{n+1}$ as required.
\end{ProofOf}

Let $\summ : \Div(E) \to E$ be the summation map.
An unexpected difficulty in the proof of Theorem~\ref{thm3}
is showing that, if $D$ is a hyperplane section for the image of $\gamma$,
then $\summ (D) = \oh$. (Comparing pull backs via 
$\mu \circ \phi_{n,n+1}$ and $\gamma \circ [n]$ only gives
that $\summ (D)$ is an $n$-torsion point.) 
We appeal to the case $N=n+1$ of the following proposition. 

\begin{Proposition}
\label{propHfixed}
Let $g \in \kbar(E)^\times$ with $\div(g) = [n]^* \oh - n^2.\oh$.
Then $g^N$ is fixed by the natural action of $H_n$ on $\L(n^2N.\oh)$.
\end{Proposition}

For the proof we need two lemmas.

\begin{Lemma}
\label{Hlem1}
Let $T \in E(\kbar)$ be a point of exact order $n$. Let
$f \in \kbar(E)^\times$ with $\div(f) = n.T-n.\oh$. If $S \in E(\kbar)$
with $2S=T$ and $nS \not= \oh$ then 
$$ \prod_{i=0}^{n-1} f(S+iT) = f(S)^n. $$
\end{Lemma}
\begin{Proof}
The rational function $P \mapsto f(P)f(T-P)$ has trivial divisor and
is therefore constant. Hence for any integer $i$,
$$ f(S+iT) f(S-iT) = f(S)^2. $$
We are immediately done in the case $n$ is odd. In the case $n$ 
is even it remains to show that $f(S) = f(S+T_2)$ where 
$T_2 = \frac{n}{2} T \in E[2](\kbar)$. Let $x \in \kbar(E)^\times$
with $\div(x) = 2.T_2 - 2.\oh$. Then comparing divisors gives
$$ \prod_{i=0}^{(n-2)/2} \tau_{-iT}^* f = c x^{n/2} $$
for some constant $c \in \kbar^\times$. Evaluating each side at 
$\pm S$ we deduce
$$ \frac{f(S+T_2)}{f(S)} = \left( \frac{x(-S)}{x(S)} \right)^{n/2} = 1. $$
\end{Proof}

\begin{Lemma}
\label{Hlem2}
Let $T \in E(\kbar)$ be a point of exact order $n$. Let
$f,g \in \kbar(E)^\times$ with $\div(f) = n.T-n.\oh$ and
$\div(g) = [n]^* \oh - n^2.\oh$. Then \\ 
(i) The pair $(f,-T)$ belongs to $H_n$ if and only if 
$\prod_{i=0}^{n-1} \tau_{iT}^* f = (-1)^{n-1}$. \\
(ii) If $f$ satisfies (i) then $\frac{g}{\tau_{-T}^* g} = f^n$. 
\end{Lemma}
\begin{Proof}
(i) By Lemma~\ref{Lemma1} the pair $(f,-T)$ belongs to $H_n$ if and only
if $(f,-T)^n = (-1)^{n-1}$. By the group law~(\ref{gplaw}) this is
equivalent to the stated condition. \\
(ii) Each side has divisor $n^2.T - n^2.\oh$. So it suffices to
check equality at $S \in E(\kbar)$ with $2S=T$ and $n S \not= \oh$.
Since $g$ has a pole of order $n^2-1$ at $\oh$, we deduce
$[-1]^* g = (-1)^{n-1} g$. So the left hand side evaluated at
$S$ is $(-1)^{n-1}$. By (i) and Lemma~\ref{Hlem1} we also get $(-1)^{n-1}$
on the right hand side.
\end{Proof}

\begin{ProofOf}{Proposition~\ref{propHfixed}}
If $(f,-T) \in H_n$ then 
$$ (f,-T) g^N = \tau_T^* (g^N/f^{nN}) = g^N $$
where for the first equality we use~(\ref{HactV*}),
and for the second equality we use Lemma~\ref{Hlem2}.
\end{ProofOf}

\section{Invariant Theory}
\label{invthy}

We recall some classical invariant theory of binary quartics and
ternary cubics, as surveyed in \cite{AKM3P}.
We then add to this theory by introducing what we call 
``covariant columns''. These are used 
in \S\ref{sec:covmat} to give formulae for $\phi_{2,3}$ and $\phi_{3,4}$.
In this section we give a complete classification of the covariant columns.
This is more that we need in~\S\ref{sec:covmat},
but serves to explain where our formulae come from.

In this section $k$ will be a field with $\charic(k) \not= 2,3$.

\subsection{Binary quartics}
\label{binq}

We study the invariants and covariants of the binary quartic
$$U(x_1,x_2) = a x_1^4 + b x_1^3 x_2 + c x_1^2 x_2^2 + d x_1 x_2^3 
+ e x_2^4. $$
For a polynomial $F \in k[x_1,x_2]$ and matrix $g \in \GL_2(k)$ we write
$$ \begin{array}{rcl}
(F \circ g)(x_1,x_2) & = & F( g_{11}x_1 + g_{12}x_2, g_{21} x_1 + g_{22} x_2)
\end{array} $$
Thus $(F \circ g)(\x) = F(g \x)$ where $\x$ 
is the column vector $(x_1,x_2)^T$. 

\begin{Definition}
\label{defcov}
A covariant $F = F(U;\x)$ of order $m$, degree $d$ and weight $p$,  
is a homogeneous
polynomial of degree $m$ in $x_1, x_2$, whose coefficients are
homogeneous polynomials of degree $d$ in the coefficients of the 
binary quartic $U$, such that
$$ F ( U \circ g ; \x) = (\det g)^p F(U; g \x) $$ 
for all $g \in \GL_2(\kbar)$. 
\end{Definition}

By considering $g$ a scalar matrix, it is clear that the 
order $m$, degree $d$ and weight $p$ of a covariant are related by
 $4d = 2p + m$.
It is well known that the ring of invariants 
(an invariant is a covariant of order~0) is
generated by $c_4$ and $c_6$ where 
$$ \begin{array}{rcl}
c_4 & = & 2^4 (12 a e - 3 b d + c^2) \\
c_6 & = & 2^5 (72ace - 27ad^2 - 27b^2e + 9bcd - 2c^3). 
\end{array} $$
Moreover, the ring of covariants is generated by $c_4$, $c_6$, 
$U$, $H$ and $J$, where
\begin{equation}
\label{Jdef}
 \begin{array}{rcl}
 H (x_1,x_2) & = & \frac{1}{3} \det \left( \frac{\partial^2 U}
{\partial x_i \partial x_j}\right)_{i,j =1,2} \end{array} $$
and $$ \begin{array}{rcl}
 J (x_1,x_2) & = & \frac{1}{12} \, 
\frac{\partial(U,H)}{\partial (x_1,x_2)} \end{array} , 
\end{equation}
subject only to the relation
\begin{equation}
\label{syz2}
 27 J^2 = -H^3 + 3 c_4 H U^2 - 2 c_6 U^3. 
\end{equation}
Since $c_4$, $c_6$, $U$, $H$ and $J$ have weights $4$, $6$, $0$, $2$ 
and $3$ we deduce
\begin{Lemma}
\label{divJ}
Every covariant of odd weight is divisible by $J$.
\end{Lemma}

We now define what we call a ``covariant column''.
\begin{Definition}
\label{defcovcol2}
A covariant column $\v = \v(U;\x)$ of order $m$, degree $d$ and weight $p$,  
is a column vector $\v=(v_1,v_2)^T$ of degree $m$ homogeneous
polynomials in $x_1, x_2$, whose coefficients are
homogeneous polynomials of degree $d$ in the coefficients of the 
binary quartic $U$, such that
$$ \v ( U \circ g ; \x) = (\det g)^p g^{-1} \v(U; g \x) $$ 
for all $g \in \GL_2(\kbar)$. 
\end{Definition}
By considering $g$ a scalar matrix, it is clear that the 
order $m$, degree $d$ and weight $p$ of a covariant column are related 
by $4d = 2p + m - 1$.
The column vector $\x$ itself is a covariant column of order $1$,
degree $0$ and weight $0$. The proof of the following lemma is
entirely straightforward, and so will be omitted.

\begin{Lemma}
(i) If $F$ is a covariant of order 
$m$, degree $d$ and weight $p$ then $\partial F = (-\frac{\partial F}
{\partial x_2}, \frac{\partial F}{\partial x_1})^T$ is a covariant
column of order $m-1$, degree $d$ and weight $p+1$, \\
(ii) If $\v_1$ and
$\v_2$ are covariant columns of orders $m_1$, $m_2$, degrees $d_1$, $d_2$
and weights $p_1$, $p_2$, then the determinant $[\v_1,\v_2]$ is
a covariant of order $m_1+m_2$, degree $d_1+d_2$ and weight $p_1+p_2-1$.
\end{Lemma}

For $F$ a covariant of order $m$ we have $[ \x, \partial F ] = m F$.

\begin{Theorem}
\label{oddeven}
(i) The covariant columns of even weight form a free 
$k[c_4,c_6,U,H]$-module
of rank $2$ with basis $\x$, $\partial J$. \\
(ii) The covariant columns of odd weight form a free 
$k[c_4,c_6,U,H]$-module
of rank $2$ with basis $\partial U$, $\partial H$. 
\end{Theorem}
\begin{Proof} (i) Since $[\x,\partial J] = 6 J \not = 0$, we can write
any covariant column as $\v = F_1 \x + F_2 \partial J$ for some
rational functions $F_1$ and $F_2$.
Then $[\v,\partial J] = 6 J F_1$ is a covariant of odd weight. It follows
by Lemma~\ref{divJ} that $F_1$ is a covariant. The same argument
shows that $F_2$ is a covariant. \\
(ii) Since we can rewrite (\ref{Jdef}) as $[\partial U,\partial H] = 12 J$,
the proof carries over exactly as in case (i).
\end{Proof}

\begin{Corollary}
The covariant columns form a module over the ring of covariants,
generated by $\x$, $\partial U$, $\partial H$, $\partial J$,
subject only to the relations 
$$ \begin{array}{rcl}
3 J \x & = & H \partial U - U \partial H \\
18 J \partial J & = & 2 (c_4 U H - c_6 U^2) \partial U 
                      + (c_4 U^2 - H^2) \partial H \\
9 J \partial U & = & (c_4 U^2- H^2) \x + 6 U \partial J \\
9 J \partial H & = & 2(c_6 U^2 - c_4 U H) \x + 6 H \partial J 
\end{array} $$
\end{Corollary}
\begin{Proof}
It only remains to describe the action of multiplication by $J$.
The first relation is obtained by applying the proof of 
Theorem~\ref{oddeven}(ii) to $\v = J \x$, and the second is
obtained by differentiating the syzygy~(\ref{syz2}).
We take linear combinations, and use the syzygy once
more to obtain the final two relations.
\end{Proof} 

\subsection{Ternary cubics}
\label{ternc}

We study the invariants and covariants of the ternary cubic
$$U(x_1,x_2,x_3) = a x_1^3 + b x_2^3 + c x_3^3 + \ldots + m x_1x_2x_3. $$
For a polynomial
$F \in k[x_1,x_2,x_3]$ and matrix $g \in \GL_3(k)$ we write
$(F \circ g)(\x) = F(g \x)$ where $\x$ is the column vector 
$(x_1,x_2,x_3)^T$. The definition of a covariant is exactly
analogous to that in the case of a binary quartic.
By considering $g$ a scalar matrix, it is clear that the 
order $m$, degree $d$ and weight $p$ of a covariant 
are related by $3d = 3p + m$.

The Hessian is a covariant of order $3$, degree $3$ and weight $2$ 
given by
$$ \begin{array}{rcl}
 H(x_1,x_2,x_3) & = & -\frac{1}{2} \det \left( \frac{\partial U}
{\partial x_i \partial x_j}\right)_{i,j =1,2,3}. \end{array} $$
There are invariants $c_4$ and $c_6$ such that
$$  H(\la U + \mu H ) =  
3 ( c_4 \la^2 \mu + 2 c_6 \la \mu^2 + c_4^2 \mu^3) U +
(\la^3 - 3 c_4 \la \mu^2 - 2 c_6  \mu^3  ) H, $$ 
and it is well known that these generate the ring of invariants.

If $Q_1(\x)$ and $Q_2(\x)$ are ternary quadrics,
with corresponding $3 \times 3$ symmetric matrices $A_1$ and $A_2$, {\em i.e.}
$$ \begin{array}{rclcrcl} Q_1(\x) & = & \frac{1}{2} \x^T A_1 \x 
&  \text{ and } & Q_2(\x) & = & \frac{1}{2} \x^T A_2 \x, \end{array} $$
then we write $\{ Q_1 , Q_2 \}$ for the $3 \times 3$ symmetric matrix
satisfying 
$$ \adj(A_1 + t A_2) = \adj(A_1) + t \{ Q_1 , Q_2 \} + t^2 \adj(A_2). $$
We can then define a symmetric matrix of quadrics
$$ \begin{array}{rcl}
M  & = & \sum_{i,j=1}^3 \left\{ \frac{\partial U}{\partial x_i} , 
\frac{\partial H}{\partial x_j} \right\} x_i x_j \end{array} $$
and covariants
\begin{equation}
\label{Jdef3}
\begin{array}{rcl} \medskip
\Theta(x_1,x_2,x_3) & = & \sum_{r,s=1}^3 M_{rs} 
\frac{\partial U}{\partial x_r}  
\frac{\partial H}{\partial x_s} \\ J(x_1,x_2,x_3) & = & \frac{1}{3} \,
\frac{\partial(U,H,\Theta)}{\partial (x_1,x_2,x_3)}. \end{array} 
\end{equation}
Again it is well known that the ring of covariants is generated
by $c_4$, $c_6$, $U$, $H$, $\Theta$ and $J$ subject only to 
a relation which reduces mod $U$ to 
\begin{equation}
\label{syz3}
 J^2 \equiv \Theta^3 - 27 c_4 \Theta H^4 - 54 c_6 H^6 \pmod{U}. 
\end{equation}

Since $c_4$, $c_6$, $U$, $H$, $\Theta$ and $J$ have 
weights $4$, $6$, $0$, $2$, $6$ and $9$ we deduce
\begin{Lemma}
\label{divJ3}
Every covariant of odd weight is divisible by $J$.
\end{Lemma}

Our definition of a covariant column is exactly analogous to that
in the case of a binary quartic. However we now also need to 
work with ``contravariant columns''.

\begin{Definition}
A covariant column, respectively contravariant column, 
$\v = \v(U;\x)$ of order $m$, degree $d$ and weight $p$,  
is a column vector $\v=(v_1,v_2,v_3)^T$ of degree $m$ homogeneous
polynomials in $x_1, x_2, x_3$, whose coefficients are
homogeneous polynomials of degree $d$ in the coefficients of the 
ternary cubic $U$, such that
$$ \v ( U \circ g ; \x) = (\det g)^p g^{-1} \v(U; g \x), $$ 
respectively
$$ \v ( U \circ g ; \x) = (\det g)^p g^T \v(U; g \x), $$ 
for all $g \in \GL_3(\kbar)$. 
\end{Definition}
By considering $g$ a scalar matrix, it is clear that the 
order $m$, degree $d$ and weight $p$ of a covariant column,
respectively contravariant column, are related 
by $3d = 3p + m - 1$, respectively $3d = 3p + m + 1$.
The column vector $\x$ itself is a covariant column of order $1$,
degree $0$ and weight $0$. The proof of the following lemma is
entirely straightforward, and so will be omitted.

\renewcommand{\labelenumi}{(\arabic{enumi})}
\begin{Lemma}
(i) If $F$ is a covariant of order 
$m$, degree $d$ and weight $p$ then $\nabla F = 
 (\frac{\partial F}{\partial x_1}, \frac{\partial F}{\partial x_2},
\frac{\partial F}{\partial x_3} )^T$ is a contravariant 
column of order $m-1$, degree $d$ and weight $p$. \\
(ii) Let $\v_1$ and $\v_2$ be covariant or contravariant columns of 
orders $m_1$, $m_2$, degrees $d_1$, $d_2$
and weights $p_1$, $p_2$. Then
\begin{enumerate}
\item 
If $\v_1$ is a covariant column, and $\v_2$ a contravariant column
then the dot product $\v_1 \cdot \v_2$ 
is a covariant of order $m_1+m_2$, degree $d_1+d_2$ and weight $p_1+p_2$. 
\item
If $\v_1$ and $\v_2$ are covariant columns then the cross 
product $\v_1 \times \v_2$ is a contravariant column of 
order $m_1+m_2$, degree $d_1+d_2$ and weight $p_1+p_2-1$. 
\item
If $\v_1$ and $\v_2$ are contravariant columns then the cross 
product $\v_1 \times \v_2$ is a covariant column of order $m_1+m_2$, 
degree $d_1+d_2$ and weight $p_1+p_2+1$.
\end{enumerate}
(iii) If $\v$ is a contravariant column of order $m$, degree $d$ and weight
$p$, then $M \v$ is covariant column of order $m+2$, degree $d+4$ and
weight $p+4$. 
\end{Lemma}

For $F$ a covariant of order $m$ we have $ \x. \nabla F  = m F$.
The determinant of the three vectors $\v_1$, $\v_2$, $\v_3$ will be denoted
$$ \begin{array}{rcccl}
[\v_1,\v_2,\v_3] & = & (\v_1 \times \v_2) \cdot \v_3 & = &
\v_1 \cdot (\v_2 \times \v_3). \end{array} $$
We define contravariant columns $\u = \nabla U$, $\h= \nabla H$, 
$\t = \nabla \Theta$ and covariant columns $\e = M\u$ and $\f = M\h$.
 
\begin{Theorem}
\label{oddeven3}
(i) The covariant columns, respectively contravariant 
columns, of even weight form a free $k[c_4,c_6,U,H,\Theta]$-module 
of rank $3$ with basis $\x$, $\e$, $\f$, respectively $\u$, $\h$, $\t$. \\
(ii) The covariant columns, respectively contravariant 
columns, of odd weight form a free $k[c_4,c_6,U,H,\Theta]$-module
of rank $3$ with basis $\u \times \h$, $\u \times \t$, $\h \times \t$,
respectively $\x \times \e$, $\x \times \f$, $\e \times \f$.
\end{Theorem}
\begin{Proof} (i)
Since $[\x,\e,\f ] = -2 J \not = 0$, we can write
any covariant column as $\v = F_1 \x + F_2 \e + F_3 \f$ for some
rational functions $F_1$, $F_2$, $F_3$.
Then $[\v,\e,\f] = -2 J F_1$ is a covariant of odd weight. It follows
by Lemma~\ref{divJ3} that $F_1$ is a covariant, and likewise for $F_2$ and 
$F_3$. The case of a contravariant column is similar, since we
can rewrite~(\ref{Jdef3}) as $[\u,\h,\t] = 3 J$. \\
(ii) 
Since $[\u \times \h, \u \times \t,\h \times \t] = [\u, \h, \t]^2 = 
9 J^2 \not= 0$,
we can write any covariant column as
$\v = F_1 (\u \times \h) + F_2 ( \u \times \t) + F_3 (\h \times \t)$ 
for some rational functions $F_1$, $F_2$, $F_3$.
Then  $\v \cdot \t = 3 J F_1$ is a covariant of odd weight. 
It follows by Lemma~\ref{divJ3} that $F_1$ is a covariant, 
and likewise for $F_2$ and $F_3$. The case of a contravariant column is 
similar.
\end{Proof}

\begin{Corollary}
(i) The covariant columns are generated as a module over the ring of covariants
by $\x$, $\e$, $\f$, $\u \times \h$, $\u \times \t$ and $\h \times \t$. \\
(ii) The contravariant columns are generated as a module over the ring of covariants 
by $\u$, $\h$, $\t$, $\x \times \e$, $\x \times \f$ and $\e \times \f$. 
\end{Corollary}

As in the case of binary quartics, there is no difficulty in finding the
relations (describing the effect of multiplication by $J$). 
Since these relations are somewhat messy, we instead record the identities
$$  \x \cdot \u  =  3 U, \quad
 \x \cdot \h  =  3 H, \quad 
 \x \cdot \t  = 6 \Theta, \quad
  \e \cdot \h  =  \f \cdot \u = \Theta $$
and 
$$ \begin{array}{rcl}
 \e \cdot \u & = & 3 (H^2 - 3 c_4 U^2) \\
 \f \cdot \h & = & 3 (3 c_4 H^2 - 8 c_6 U H + 3 c_4^2 U^2) \\
 \e \cdot \t & = & 12 (3 c_4 H^3 - c_4 U \Theta - 12 c_6 U H^2 + 9 c_4^2 U^2 H)  \\
 \f \cdot \t & = & 12 (c_4 H \Theta - 3 c_6 H^3  - 3 c_6 U \Theta 
    + 9 c_4^2 U H^2 - 3 c_4 c_6 U^2 H - 9 c_4^3 U^3 )
\end{array} $$
from which the relations may be recovered by following the 
proof of Theorem~\ref{oddeven3}. 

\section{Covariant matrices}
\label{sec:covmat}

We use the invariant theory of binary quartics and ternary cubics
to given an alternative proof of Theorem~\ref{thm3} in 
the cases $n=2,3$.  

\subsection{The case $n=2$}
A binary quartic $U \in k[x_1,x_2]$ is non-singular if
its discriminant $\Delta = (c_4^3-c_6^2)/1728$ is non-zero.
We write $\PP(1,1,2)$ for the weighted projective space
where the co-ordinates $x_1$, $x_2$, $y$ are assigned 
degrees $1$, $1$, $2$.

\begin{Proposition}
\label{prop:covmat2}
Let $U \in k[x_1,x_2]$ be a non-singular binary quartic with invariants
$c_4, c_6 \in k$ and covariants $H, J \in k[x_1,x_2]$. Then \\
(i) The equation $y^2 = U(x_1,x_2)$ defines a smooth curve of genus
one $C_2 \subset \PP(1,1,2)$. \\
(ii) The Jacobian $E$ of $C_2$ has Weierstrass equation
$$ Y^2 Z = X^3 - 27 c_4 X Z^2 - 54 c_6 Z^3. $$
(iii) The $2$-covering map $\pi: C_2 \to E$ is given by
$$ (Z:X:Y) = (y U(x_1,x_2) : -3 y H(x_1,x_2) : 27 J(x_1,x_2) ). $$
(iv) Let $\phi_3 : E \to \PP^2$ be the natural inclusion.
Then the  Brauer-Severi diagrams $\phi_2 : C_2 \to \PP^1$ $;$ 
$(x_1:x_2:y) \mapsto (x_1:x_2)$ and
$ \phi_{2,3}:  C_2  \to  \PP(\Mat_{2,3})$ $;$ 
$(x_1:x_2:y) \mapsto A_{2,3}$, where
$$ A_{2,3} = \begin{pmatrix} \smallskip
-9 \frac{\partial H}{\partial x_2} &
 -3 \frac{\partial U}{\partial x_2} & x_1 y \\
\,\,\,\, 9 \frac{\partial H}{\partial x_1} & 
\,\, \,\, 3 \frac{\partial U}{\partial x_1} & x_2 y \\
\end{pmatrix}, $$ 
satisfy the conclusions of Theorem~\ref{thm3}.
\end{Proposition}

\begin{Proof}
Statements (i)--(iii) are well known: see \cite{AKM3P}, \cite{Weil}. \\
(iv) Let $D$ be the hyperplane section for $\phi_2 : C \to \PP^1$.
Then $\L(3D)$ has basis $x_1^3$, $x_1^2 x_2$, $x_1 x_2^2$, $x_2^3$, $x_1y$,
$x_2y$. We write the entries of $A_{2,3}$ as linear combinations of these 
basis elements, and arrange the coefficients in a $6 \times 6$ matrix.
The determinant of this matrix is $2^2 3^8 \Delta$. Hence 
$\phi_{2,3}: C_2 \to \PP(\Mat_{2,3})$ is an embedding by the complete
linear system $|3D|$.

In the notation of \S\ref{binq} we have $A_{2,3} = (9 \partial H, 3 \partial U,
 y \x )$. The $2 \times 2$ minors of this matrix are 
$$ \begin{pmatrix} -3 y [\x ,\partial U] \\ 9 y [\x ,\partial H] \\
-27 [\partial U, \partial H] \end{pmatrix} = 12 \begin{pmatrix}
-y U \\ 3 y H \\ -27 J \end{pmatrix} $$
So the final statement of Theorem~\ref{thm3} is immediate from (iii).

It remains to show that $S \in E[2](\kbar)$, respectively $T \in E[3](\kbar)$,
acts on the image of $\phi_{2,3}$ as left multiplication by $\chi_2(S)$,
respectively right multiplication by $\chi_3(T)$.

The statement for $S \in E[2](\kbar)$ follows formally from the
covariance of the columns of $A_{2,3}$. Indeed, writing
$A_{2,3} = (\a_1, \a_2, y \a_3 )$, where the $\a_i$ are covariant columns,
Definition~\ref{defcovcol2} gives
$$ \a_i ( U \circ g ; \x) = g^{-1} \a_i (U, g \x) $$
for all $g \in \SL_2(\kbar)$. So if $g \in \SL_2(\kbar)$ is a lift
of $\chi_2(S)$ then $U \circ g = U$ and 
$$ g \, \phi_{2,3}(\x: y) = \phi_{2,3} (g \, \x : y) $$
as required.

We deduce the statement for $T \in E[3](\kbar)$ from the parts
of Theorem~\ref{thm3} already established. To this end, let $M_S$ 
and $M_T$ be endomorphisms of $\Mat_{2,3}(\kbar)$ lifting
$\chi_{2,3}(S)$ and $\chi_{2,3}(T)$ respectively. We have
shown that $M_S$ is left multiplication by a $2 \times 2$ matrix.
Accordingly we view $\Mat_{2,3}(\kbar)$ as an $H_2$-module
via left multiplication. (In this proof we write $H_n$ 
as a shorthand for $H_n(\kbar)$.)
Since 2 and 3 are coprime, the commutator of $M_S$ and $M_T$ 
is trivial. Hence $M_T$ is an endomorphism of $\Mat_{2,3}(\kbar)$ 
as an $H_2$-module. Since the standard representation 
$V_2$ of $H_2$ is irreducible, it follows by Schur's lemma that 
$M_T$ belongs to 
$$ \End_{H_2} (\Mat_{2,3}(\kbar)) \isom 
 \End_{H_2} ( V_2 \oplus V_2 \oplus V_2 ) \isom \Mat_3(\kbar). $$
Thus $M_T$ is right multiplication by a $3 \times 3$ matrix.
Let $\psi(T)$ be the image of this matrix in  $\PGL_3(\kbar)$. 
It remains to show that the characters 
$\chi_3: E[3] \to \PGL_3$ and $\psi : E[3] \to \PGL_3$ are equal.
Recalling that $\phi_3 \circ \pi = \mu \circ \phi_{2,3}$, we take
$P \in C_2(\kbar)$ and compute
$$ \begin{array}{lrcl}
& \phi_3 (\pi(P+T)) & = & \mu (\phi_{2,3}(P+T)) \\
\implies & \phi_3 (\pi (P) + 2 T) & = & \mu (\phi_{2,3}(P) \psi(T)) \\
\implies & \chi_3(T)^{-1} \phi_3 (\pi (P))
& = & \psi(T)^{-1}  \mu (\phi_{2,3}(P)).
\end{array} $$
Hence $\chi_3 = \psi$ as required.
\end{Proof}

\subsection{The case $n=3$}

A ternary cubic $U \in k[x_1,x_2,x_3]$ is non-singular if
its discriminant $\Delta = (c_4^3-c_6^2)/1728$ is non-zero.
The covariant columns $\e$, $\f$ and contravariant columns
$\u$, $\h$ were defined in \S\ref{ternc}. 

\begin{Proposition}
\label{prop:covmat3}
Let $U \in k[x_1,x_2,x_3]$ be a non-singular ternary cubic with invariants
$c_4, c_6 \in k$ and covariants $H, \Theta, J \in k[x_1,x_2,x_3]$. Then \\
(i) The equation $U(x_1,x_2,x_3)=0$ defines a smooth curve of genus
one $C_3 \subset \PP^2$. \\
(ii) The Jacobian $E$ of $C_3$ has Weierstrass equation
$$ Y^2 Z = X^3 - 27 c_4 X Z^2 - 54 c_6 Z^3. $$
(iii) The $3$-covering map $\pi: C_3 \to E$ is given by
$$ (Z:X:Y) = (H^3 : \Theta H : J ). $$
(iv) Let $\phi_4 : E \to \PP^3$ $;$ $(Z:X:Y) \mapsto (Z^2:XZ:YZ:X^2)$. 
Then the Brauer-Severi diagrams $\phi_3 : C_3 \to \PP^2$ and 
$ \phi_{3,4}:  C_3  \to  \PP(\Mat_{3,4})$ $;$ 
$(x_1:x_2:x_3) \mapsto A_{3,4}$, where
$$ A_{3,4} = \begin{pmatrix}
-3 \f + 9 c_4 H \x & \e & \frac{2}{3} (\u \times \h) & -\frac{1}{3} H \x 
\end{pmatrix}, $$
satisfy the conclusions of Theorem~\ref{thm3}.
\end{Proposition}

\begin{Proof}
Statements (i)--(iii) are well known: see \cite{AKM3P}. \\
(iv) By Theorem~\ref{oddeven3} the covariant columns of
order 4 form a free $k[c_4,c_6]$-module with basis
\begin{equation}
\label{order4basis}
 U \x, \quad H \x, \quad \e, \quad \f, \quad \u \times \h.
\end{equation}
The entries of these columns give us 15 ternary quartics.
We arrange the coefficients of these quartics in a 
$15 \times 15$ matrix, and find that the determinant is
$2^{42}3^{12} \Delta^5$. (The calculation is made easier if we first 
put $U$ in Hesse normal form:
$$ U(x_1,x_2,x_3) = a(x_1^3 + x_2^3 + x_3^3) + b x_1 x_2 x_3.) $$

Let $D$ the hyperplane section for $\phi_3 : C_3 \to \PP^2$.
Since the only ternary quartics vanishing on $C_3$ are the 
entries of $U \x$,  the above calculation shows that 
$\phi_{3,4}: C_3 \to \PP(\Mat_{3,4})$ is an embedding by the 
complete linear system $|4D|$.

A direct calculation (carried out for $U$ in Hesse normal form)
shows that the $3 \times 3$ minors of $A_{3,4}$ are
$$ \begin{array}{rcl} 
\mu_1 & = & 2 H^4 - 6 c_4 U^2 H^2 - \frac{2}{3} U H \Theta \\
\mu_2 & = & 2 \Theta H^2 - 18 c_4^2 U^3 H - 18 c_4 U H^3 + 48 c_6 U^2 H^2 \\
\mu_3 & = & 2 J H \\
\mu_4 & = & 2 \Theta^2 + 162 c_4^3 U^4 - 54 c_4^2 U^2 H^2
   - 432 c_4 c_6 U^3 H \qquad \\
&  & \multicolumn{1}{r}{ - \, 18 c_4 U H \Theta + 144 c_6 U H^3.}
\end{array} $$
The final statement of Theorem~\ref{thm3} follows since
by (iii) the composition $\phi_4 \circ \pi$ is given by
$$(x_1:x_2:x_3) \mapsto (H^4 : \Theta H^2 : J H : \Theta^2).$$ 

The remainder of the proof now carries through exactly as
in the case $n=2$.
\end{Proof}

\section{Computations}
\label{implement}

Let $E$ be an elliptic curve over $\Q$. We use $6$-descent and
$12$-descent to assist in the search for generators of $E(\Q)$
of large height. The method is of greatest interest when $E(\Q)$ 
has rank at least 2, or $E$ has large conductor, {\em i.e.}
in those cases where we cannot use Heegner points.
In this section we give some details of our implementation in the
computer algebra system Magma \cite{Magma}. Further remarks accompany
the numerical examples in \S\ref{sec:num}.

\subsection{The method in outline}
\label{sec:outline}

We begin by using the existing functions in Magma 
to compute $n$-coverings for $n=2,3,4$.

\begin{itemize}
\item The Magma function {\tt TwoDescent}, takes
as input a Weierstrass equation for $E$, and returns a list of
$2^s-1$ binary quartics representing the non-zero elements of 
the $2$-Selmer group \\ $S^{(2)}(E/\Q) \isom (\Z/2\Z)^s$. 
\item The Magma function {\tt ThreeDescent}, written by Stoll, Donnelly and
the author, takes as input a Weierstrass equation for $E$, 
and returns a list of
$(3^t-1)/2$ ternary cubics representing the non-zero elements of 
the $3$-Selmer group $S^{(3)}(E/\Q) \isom (\Z/3\Z)^t$. 
\item The Magma function {\tt FourDescent}, written by Womack and Watkins,
takes as input a binary quartic representing a non-zero element 
$\alpha \in S^{(2)}(E/\Q)$, and returns a list of pairs of quadrics 
in four variables, representing the elements of the $4$-Selmer group 
in the fibre of
$S^{(4)}(E/\Q) \to S^{(2)}(E/\Q)$ above $\alpha$.
\end{itemize}

Each element of the $n$-Selmer group is now 
represented by (equations for) a Brauer-Severi diagram $[C_n \to \PP^{n-1}]$.
The Selmer group elements may equally be viewed as $n$-coverings, where
the covering maps $\pi: C_n \to E$ are computed using the classical
formulae surveyed in \cite{AKM3P}. (For $n=2$, $3$ we recalled these formulae
in Propositions~\ref{prop:covmat2} 
and~\ref{prop:covmat3}.) Replacing the covering map $\pi$ by
$[-1] \circ \pi$ corresponds to taking the inverse in the Selmer group.
So in the cases $n=3,4$ each ternary cubic, respectively pair 
of quadrics, represents both a Selmer group element and its inverse. 

\begin{Definition} 
Let $[C_n \to \PP^{n-1}]$ be a Brauer-Severi diagram
with hyperplane section $D$. A point $P \in C_n(\kbar)$ is a flex 
if $n.P \sim D$. 
\end{Definition}
The flex points of a 2-covering are the roots of the binary quartic.
In the cases $n=3,4$ the flex points (also known as points of inflection,
or hyperosculating points) are the intersections with $H=0$,
respectively $J=0$, where $H$ is the Hessian of a ternary cubic, 
and $J$ is the covariant defined in \cite[\S3.3]{AKM3P}. 

We recall that if $[C_n \to \PP^{n-1}]$ is a Brauer-Severi diagram,
then the morphism $C_n \to \PP^{n-1}$ is that determined by a complete
linear system of degree $n$. So if $n \ge 3$ then $C_n \to \PP^{n-1}$
is an embedding. We identify $C_n$ with its image, which is called
a genus one normal curve of degree~$n$. It is well known that if
$n \ge 4$ then the homogeneous ideal $I(C_n)$ is generated by a 
vector space of quadrics of dimension $n(n-3)/2$.

The details of $6$-descent are as follows. We start with a $2$-covering
$C_2 = \{ y^2 = U_2(x_1,x_2) \}$ and a $3$-covering 
$C_3 = \{ U_3(x_1,x_2,x_3) = 0 \}$, each defined over $\Q$. 
Since these are coverings of the same elliptic curve $E$,
we may assume that $U_2$ and $U_3$ have the same invariants 
$c_4$ and $c_6$. Then $E$ has Weierstrass equation
$$ y^2 = x^3 - 27 c_4 x - 54 c_6. $$
We compute a flex point on $C_3$ with co-ordinates in a number field, 
$L$ say. Typically $[L: \Q]=9$. Then Algorithm~\ref{flexmat3} 
finds a matrix $g \in \GL_3(L)$ with
$$ (U_3 \circ g) (z,x,y) 
= \lambda ( y^2 z - x^3 + 27 c_4 x z^2 + 54 c_6 z^3)$$
for some $\lambda \in L^\times$. (In the notation of \S\ref{sec:comptwists} 
we have $g = B_3^{-1}$.) 

Next we let $\phi_{2,3} : C_2 \to \PP(\Mat_{2,3})$ 
be the embedding defined in Proposition~\ref{prop:covmat2}(iv).
The image is a genus one normal curve of degree~$6$. We use linear
algebra to compute a basis $Q_1, \ldots, Q_9$ for the space
of quadrics vanishing on this curve. Writing these as
polynomials in variables $X_{ij}$ for $1 \le i \le 2$ and $1 \le j \le 3$,
we make the substitution
$$  \begin{pmatrix} X_{11} & X_{12} & X_{13} \\
 X_{21} & X_{22} & X_{23} \end{pmatrix} = 
 \begin{pmatrix} x_{11} & x_{12} & x_{13} \\
 x_{21} & x_{22} & x_{23} \end{pmatrix} g. $$
The new quadrics have coefficients in $L$, but the vector space
they span has a basis with coefficients in $\Q$. We compute
an LLL-reduced basis for the intersection of this space with
$\Z[x_{11}, \ldots, x_{23}]$. These are now the equations for a 
6-covering $C_6 \subset \PP(\Mat_{2,3})$. Moreover, by
Theorem~\ref{thm3}, the covering map $C_6 \to C_3$ is defined 
by the $2 \times 2$ minors, {\em i.e.} 
$$ \begin{pmatrix} x_{11} & x_{12} & x_{13} \\
 x_{21} & x_{22} & x_{23} \end{pmatrix} 
\mapsto ( x_{12} x_{23} - x_{22} x_{13} : x_{13} x_{21} - x_{23} x_{11}
: x_{11} x_{22} - x_{21} x_{12} ). $$

The details of 12-descent are similar.
We start with a $3$-covering 
$C_3 = \{ U_3(x_1,x_2,x_3) = 0 \}$ and a $4$-covering
$C_4 = \{ Q_1 = Q_2 = 0 \}$, each defined over $\Q$, and with
the same invariants $c_4$ and $c_6$. We compute a flex point on $C_4$ 
with co-ordinates in a number field, $L$ say. Typically $[L: \Q]=16$. 
Then Algorithm~\ref{flexmat4} finds a matrix $g \in \GL_4(L)$ with
$$ \langle Q_1 \circ g, Q_2 \circ g \rangle = \langle x_1 x_4 - x_2^2,
x_2 x_4 - x_3^2 - 27 c_4 x_1 x_2 -54 c_6 x_1^2 \rangle.  $$

Next we let $\phi_{3,4} : C_3  \to \PP(\Mat_{3,4})$ be the embedding 
defined in Proposition~\ref{prop:covmat3}(iv).
The image is a genus one normal curve of degree~$12$. We use linear
algebra to compute a basis $Q_1, \ldots, Q_{54}$ for the space
of quadrics vanishing on this curve. As in the case of $6$-descent,
we then twist by $g \in \GL_4(L)$ to obtain equations for a 
12-covering $C_{12} \subset \PP(\Mat_{3,4})$. 
Moreover, by Theorem~\ref{thm3},
the covering map $C_{12} \to C_4$ is defined by the $3 \times 3$ minors.

Unlike the case of 6-descent, we can combine a 3-covering and a 4-covering
to give a 12-covering in two essentially different ways. This is
because each of $C_3$ and $C_4$ represents both a Selmer group 
element and its inverse. It is important that we compute both 
12-coverings, since in the case they are soluble, 
their rational points will cover $\pm P + 12 E(\Q)$
and $\pm 5P + 12 E(\Q)$ respectively. 
In practice the second 12-covering
is obtained by switching the sign in the third column of the matrix
defining $\phi_{3,4}$. 

It remains to search for rational points on $C_6$ and $C_{12}$.
We use the $p$-adic point searching method due independently to
Elkies and Heath-Brown, as implemented by Watkins in 
the Magma function {\tt PointSearch}. Descriptions may be found in 
\cite{Watkins} and \cite[\S2.9]{Womack}. (Elkies' original 
paper~\cite{Elk00} only considers real approximations.)
The method first chooses an auxiliary
prime $p$, whose size depends on the height bound set for the
search. The points on the reduction of $C$ mod $p$ are then
enumerated, and for each such point $P_0$ a lattice method variant
of Hensel's lemma is used to search for rational points on $C$
with reduction $P_0$. A variant of the method uses two primes.
The method works particularly well for curves of high codimension 
as considered here. 

Finally, our search for points is significantly improved if we ``minimise''
our equations for $C_6$ and $C_{12}$ before running {\tt PointSearch}.
We give details in \S\ref{sec:min}.

\subsection{Computing flex matrices}
\label{compfm}

Let $[C \to \PP^{n-1}]$ be a Brauer-Severi diagram. 
To compute a flex matrix for $C$, as defined in \S\ref{sec:comptwists}, we 
first find a flex point $P$ on $C$. 
We then follow an inductive procedure, based on the idea of 
projecting away from $P$. This method is a by-product 
of the standard procedures for putting an elliptic curve
in Weierstrass form, as described in \cite[\S8]{CaL}. 
We therefore simply list the algorithms used. Notice that we do not 
use the general Riemann-Roch machinery implemented in Magma,
as this would be unnecessarily slow in our applications.

\begin{Algorithm}
\label{flexmat2}
Let $U \in k[x_1,x_2]$ be a non-singular binary quartic with
invariants $c_4$ and $c_6$. Given $(\alpha : \beta) \in \PP^1(k)$
with $U(\alpha,\beta)=0$ we compute $g \in \GL_2(k)$ with last
column $(\alpha,\beta)^T$ satisfying
$$ \begin{array}{rcl} (U \circ g)(z,x) & = & \frac{1}{36} (\det g)^2 
 (x^3 z - 27 c_4 x z^3 - 54 c_6 z^4). \end{array} $$
\begin{enumerate}
\item Choose any $g_1 \in \GL_2(k)$ with last column $(\alpha,\beta)^T$.
\item Compute $(U \circ g_1)(z,x) = (\det g_1)^2 
( b x^3 z + c x^2 z^2 + \ldots)$ and put $g_2 = \begin{pmatrix} 36b & 0 \\
-12c & 1 \end{pmatrix}$.
\item Return $g_1g_2$.
\end{enumerate}
\end{Algorithm}

\begin{Algorithm}
\label{flexmat3}
Let $U \in k[x_1,x_2,x_3]$ be a non-singular ternary cubic with
invariants $c_4$ and $c_6$. Given $(\alpha : \beta:\gamma) \in \PP^2(k)$
a flex point on the curve $\{U=0\}$ we compute $g \in \GL_3(k)$ with last
column $(\alpha,\beta,\gamma)^T$ satisfying
$$  \begin{array}{rcl} (U \circ g)(z,x,y) & = & \frac{1}{6} (\det g) 
 ( y^2 z - x^3 + 27 c_4 x z^2 + 54 c_6 z^3 ). \end{array} $$
\begin{enumerate}
\item Choose any $g_1 \in \GL_3(k)$ with last column $(\alpha,\beta,\gamma)^T$.
\item Write $(U \circ g_1)(z,x,y) = (\det g_1) 
( f_1(z,x) y^2 + f_2(z,x)y + f_3(z,x) )$ 
and let $\alpha, \beta \in k$ with $f_1(z,x) = \beta z - \alpha x$. 
Then run Algorithm~\ref{flexmat2}
on $\frac{1}{4} f_2^2 - f_1 f_3$ to obtain $g \in \GL_2(k)$.
\item Compute $p,q \in k$ with
$$ \begin{array}{rcl}
(f_1 \circ g)(z,x) & = & (\det g) z \\
(f_2 \circ g)(z,x) & = & (\det g) (p z + q x) z \\
\end{array} $$
and put $g_2 = \begin{pmatrix} 6g & \begin{array}{c} 0 \\ 0 \end{array} \\
\begin{array}{cc} -3p & -3q \end{array} & 1 \end{pmatrix}$.
\item Return $g_1g_2$.
\end{enumerate}
\end{Algorithm}

In the case $n=4$ the invariants are again described in \cite{AKM3P}.
We label them $c_4$ and $c_6$, with scalings as specified in 
\cite{g1inv}.

\begin{Algorithm}
\label{flexmat4}
Let $Q_1,Q_2 \in k[x_1,x_2,x_3,x_4]$ be a pair of quadrics, 
with invariants $c_4$ and $c_6$. We suppose that  $\{Q_1=Q_2=0\}$ 
is a smooth curve of genus one. 
Given $(\alpha : \beta: \gamma : \delta) \in \PP^3(k)$,
a flex point on this curve, we compute $g \in \GL_4(k)$ 
with last column $(\alpha,\beta,\gamma,\delta)^T$ satisfying
$$ \langle Q_1 \circ g, Q_2 \circ g \rangle = \langle x_1 x_4 - x_2^2,
x_2 x_4 - x_3^2 - 27 c_4 x_1 x_2 -54 c_6 x_1^2 \rangle.  $$
\begin{enumerate}
\item Choose any $g_1 \in \GL_4(k)$ with last 
column $(\alpha,\beta,\gamma,\delta)^T$.
\item Write 
$$ \begin{array}{rcl}
Q_1 \circ g_1 & = & \ell_1(x_1,x_2,x_3) x_4 + q_1(x_1,x_2,x_3) \\
Q_2 \circ g_1 & = & \ell_2(x_1,x_2,x_3) x_4 + q_2(x_1,x_2,x_3) 
\end{array} $$
and let 
$$     \alpha = \left| \begin{matrix} \ell_{12} & \ell_{13} 
\\ \ell_{22} & \ell_{23} \end{matrix} \right|, 
\quad   \beta = \left| \begin{matrix} \ell_{13} & \ell_{11} 
\\ \ell_{23} & \ell_{21} \end{matrix} \right|, 
 \quad \gamma = \left| \begin{matrix} \ell_{11} & \ell_{12} 
\\ \ell_{21} & \ell_{22} \end{matrix} \right|,  $$
where $\ell_i = \sum \ell_{ij} x_j$. Then run Algorithm~\ref{flexmat3}
on $(\det g_1)^{-1}(\ell_2 q_1 - \ell_1 q_2)$ to obtain $g \in \GL_3(k)$.
\item Replace $Q_1$ and $Q_2$ by linear combinations (and update
the $\ell_i$ and $q_i$ of Step 2) so that $\ell_i \circ g = x_i$ for
$i=1,2$.
\item Compute $a,b,c \in k$ with
$$ \begin{array}{rcl} \smallskip
\qquad q_1 \circ g & = & \frac{1}{6}(x_1(a x_1 + b  x_2 + c x_3) - x_2^2)  \\
q_2 \circ g & = & \frac{1}{6}(x_2(a x_1 + b  x_2 + c x_3) 
- x_3^2  - 27 c_4 x_1 x_2 - 54 c_6 x_1^2) \\
\end{array} $$
and put $g_2 = \begin{pmatrix} 6g & \begin{array}{c} 0 \\ 0 \\ 0 \end{array} \\
\begin{array}{ccc} -a & -b & -c \end{array} & 1 \end{pmatrix}$.
\item Return $g_1g_2$.
\end{enumerate}
\end{Algorithm}

\subsection{Minimisation}
\label{sec:min}

If an $n$-covering $\pi : C \to E$ is to be useful in the
search for rational points on $E$, not only must we find 
explicit equations for $C \subset \PP^{n-1}$, 
but we must also find a change of
co-ordinates on $\PP^{n-1}$ so that these equations 
have reasonably small coefficients.
The task naturally falls into two parts which, following terminology
introduced by Cremona, we call minimisation and reduction. 

Minimisation is the task of removing as many prime factors as
possible from a suitably defined discriminant. The most familiar
example is that of minimising a Weierstrass equation. 
By reduction we mean the use of unimodular transformations
to further decrease the size of the coefficients. The basic example
is reduction of binary quadratic forms, or 
more generally lattice reduction.
Thus minimisation is concerned with the finite places, and
reduction with the infinite places. 
The need to perform reduction
is our main reason for working over the rationals 
(instead of a more general number field).

The minimisation and reduction of 2-coverings has been studied
in \cite{BSD}, \cite{CreRed}, and \cite{CS}.
The generalisations to $3$-coverings and $4$-coverings 
are described in \cite{paperIV} and \cite{Womack}. 
These algorithms have been implemented in Magma, and are
called by the functions {\tt TwoDescent}, {\tt ThreeDescent} and 
{\tt FourDescent}. Hence in \S\ref{sec:outline} we start with
an $n$-covering and an $(n+1)$-covering both of which are already
minimised and reduced. So it would not be unreasonable to hope 
that the $n(n+1)$-covering computed from them will automatically 
be minimised and reduced. Numerical examples
suggest that this is true for reduction, but not for minimisation.

The following is a description of our current ad hoc approach to 
the minimisation of $n$-coverings for $n>5$. 
Although this method works reasonably well in practice, 
there remains considerable room for both theoretical and 
practical improvements.
\enlargethispage{0.1cm}

Let $C \subset \PP^{n-1}$ be a genus one normal curve, defined
over $\Q$, and of degree $n \ge 4$. 
We recall that the homogeneous ideal $I(C)$ in $\Q[x_1, \ldots, x_n]$
is generated by a vector space of quadrics of dimension
$N = n(n-3)/2$. Then for $A \in \GL_N(\Q)$ and $B \in \GL_n(\Q)$ 
we define
$$ \begin{array}{rcl} 
[A,B](q_1, \ldots, q_N) & = & ( \sum_{i=1}^N a_{1i} q_i \circ B^T, 
\ldots , \sum_{i=1}^N a_{Ni} q_i \circ B^T ). 
\end{array} $$

\begin{Definition}
\label{def:min}
(i) An integral model $(q_1, \ldots, q_{N})$ for $C$ is a tuple 
of quadrics in $\Z[x_1, \ldots, x_n]$ generating $I(C)$. \\
(ii) An integral model $(q_1, \ldots, q_{N})$ is minimal at a prime $p$
if whenever $A \in \GL_N(\Q)$ and $B \in \GL_n(\Q)$ with
$[A,B](q_1, \ldots, q_N)$ integral, then
$$ \ord_p(\det A) + ( n - 3 ) \ord_p(\det B) \ge 0. $$
(iii) An integral model is minimal if it is minimal at all primes $p$.
\end{Definition}

It is not even clear from our definition that minimal models exist,
although in the cases $n=4$, $5$ this can be proved using the
invariants defined in~\cite{AKM3P} and~\cite{g1inv}.
The condition in Definition~\ref{def:min}(ii) is motivated by 
considering what happens when $A$ and $B$ are scalar matrices.

We attempt to minimise at $p$ as follows. Let $I_p$ be the ideal
in $\F_p[x_1, \ldots, x_n]$ generated by the reductions 
of $q_1,\ldots,q_N$ mod $p$. We compute the $\F_p$-vector space
$V_p$ of linear forms in the radical of $I_p$. By a unimodular change
of co-ordinates we may suppose that $V_p = \langle x_1, \ldots, x_d \rangle$
for some $0 \le d < n$. Then we put
$$ B = \begin{pmatrix} p I_d & 0 \\ 0 & I_{n-d} \end{pmatrix} $$
and compute the index $p^m$ of the lattice spanned by the $q_i \circ B$
in its saturation in $\Z[x_1,\ldots,x_n]$. We call the integer
$m - (n-3)d$ the gain. If the gain is positive then 
we switch to the new quadrics and start over again. 
Otherwise we stick with the old quadrics. 
Sometimes it is worth trying other choices for the 
vector space $V_p$, for example the space of linear forms in one 
of the minimal primes containing $I_p$. There is no guarantee that 
these methods will produce a $p$-minimal model (and in general
they do not).

The end result of our attempts at minimisation is a 
change of co-ordinates on $\PP^{n-1}$. We run the LLL algorithm 
on the rows (or columns depending on conventions) of the change
of basis matrix, before applying it to the original quadrics.
This is to ensure that we do not throw away the fact our quadrics 
are already (close to being) reduced. 

\section{Numerical Examples}
\label{sec:num}

\subsection{An example of 6-descent}

The following problem falls into the class of problems discussed 
on pages 480-481 of \cite{Dickson}.
\begin{quotation}
Given an integer $N>2$, decide whether there are non-zero integers
$x$ and $y$ such that both $x^2+Nxy+y^2$ and $x^2-Nxy+y^2$ are squares.
\end{quotation}
Elementary manipulations show that the problem is equivalent to deciding
whether the elliptic curve
$$   E_N : \quad  y^2  =   x ( x + ( N + 2)^2 ) ( x + ( N - 2)^2 )  $$
has positive rank\footnote{The trivial solutions correspond to a subgroup
$T \subset E_N(\Q)$ with $T \isom \Z/2\Z \times \Z/4\Z$. Since $\pm(N^2-4)
\notin (\Q^*)^2$, the image of $T$ under the $2$-descent map
$E_N(\Q) \to \Q^*/(\Q^*)^2 \times \Q^*/(\Q^*)^2$ has order $4$. 
It follows that $E_N(\Q)_{\tors} = T$.}. MacLeod and Rathbun \cite{MacLeod} 
have undertaken to find a solution for $x$ and $y$ 
(where one exists) for all $N < 1000$.
The one case to elude them (as of November 2006) was $N=809$,
for which the rank is 1 and the generator
is predicted\footnote{This estimate comes from the Birch--Swinnerton-Dyer 
conjecture, assuming that the Tate-Shafarevich group is trivial.}
to have height $617.88$. On the 2-isogenous curve
$$ E'_N : \quad  y^2 = x^3 + 2 (N^2+12N+4) x^2 + (N-2)^4 x $$
the predicted height is half this value, yet still beyond the range that
can be found using 4-descent. The conductor of $E = E_{809}'$ 
is sufficiently large that a Heegner point calculation 
ran into difficulties (and for this reason the curve was reported to
Magma as a bug).

We find a point of infinite order on $E$ using 6-descent.
The existing Magma functions for $2$-descent and
$3$-descent give us a 2-covering
$$ C_2 = \left\{ y^2 =  138546 x_1^4 
  + 225978 x_1^3 x_2 + 435649 x_1^2 x_2^2
   + 3884 x_1 x_2^3 + 183499 x_2^4 \right\} $$
and a 3-covering 
$$ C_3 = \left\{ \begin{array}{rcl} 
54 x^3 - 84 y^3 - 258 z^3  + 144 x^2 y + 87 x^2 z  - 350 x y^2 \quad \\
+ \, 71 y^2 z - 1656 x z^2 - 986 y z^2  - 388 x y z 
& = & 0 \end{array} \right\}. $$
To compute this 3-covering we had to find the class group and units for 
a number field of degree 8. This is by far the most time consuming part 
of the $6$-descent calculation, taking a couple of hours, 
as compared to at most a couple of minutes for each of the other steps. 

Following the method described in \S\ref{sec:outline} we compute
9 quadrics defining a 6-covering $C_6 \subset \PP(\Mat_{2,3})$.
These are quadrics in 6 variables,
labelled $x_{ij}$ for $1 \le i \le 2$ and $1 \le j \le 3$.
The coefficients are reasonably small integers, the largest in absolute
value being 142.

Minimising (at the primes $2$, $3$, $809$ and $811$), 
as described in \S\ref{sec:min}, 
suggests making the substitution
$$ \begin{pmatrix} x_{11} \\ x_{12} \\ x_{13} \\ x_{21} \\ x_{22} \\ x_{23}
\end{pmatrix} = 
\begin{pmatrix}
     -70 & 455 & 293 & 700 & -63 & 437 \\
     104 & 417 & -363 & 290 & 745 & -579 \\
     -268 & -89 & 205 & -60 & 1223 & 817 \\
     -320 & -335 & -839 & 386 & 147 & -311 \\
     -318 & 411 & -123 & -696 & -405 & 561 \\
     284 & 59 & -523 & 226 & -15 & 1973 
\end{pmatrix}
\begin{pmatrix} x_1 \\ x_2 \\ x_3 \\ x_4 \\ x_5 \\ x_6
\end{pmatrix}. $$
This decreases the size of the coefficients.
More importantly, but as we will only see in hindsight, 
it also reduces the (naive) height of the point we are looking for. 
The new quadrics are
\tiny
$$
\begin{array}{rcl}
q_1 & = & x_{1}x_{4}+2x_{2}x_{3}+2x_{2}x_{4}-2x_{2}x_{5}-2x_{2}x_{6}+2x_{3}^2-2x_{3}x_{4}-6x_{3}x_{6}-x_{4}x_{5}+4x_{4}x_{6} \\ 
 & & \quad -2x_{5}^2+5x_{5}x_{6}+x_{6}^2 \\ 
q_2 & = & x_{1}x_{2}+x_{1}x_{3}-2x_{1}x_{4}+3x_{1}x_{5}-4x_{1}x_{6}-x_{2}^2+2x_{2}x_{5}-7x_{2}x_{6}+x_{3}^2-x_{3}x_{4}+4x_{3}x_{5} \\ 
 & & \quad -2x_{3}x_{6}+5x_{4}^2-3x_{4}x_{5}-8x_{4}x_{6}+x_{5}^2+2x_{5}x_{6}-6x_{6}^2 \\ 
q_3 & = & x_{1}^2+x_{1}x_{3}+3x_{1}x_{4}+5x_{1}x_{5}-2x_{1}x_{6}-3x_{2}^2+5x_{2}x_{3}+x_{2}x_{4}-4x_{2}x_{5}-x_{2}x_{6}+2x_{3}^2 \\ 
 & & \quad +4x_{3}x_{4}+7x_{3}x_{5}-5x_{3}x_{6}+3x_{4}^2-5x_{4}x_{6}-2x_{5}x_{6}-x_{6}^2 \\ 
q_4 & = & x_{1}^2-3x_{1}x_{2}+3x_{1}x_{4}+x_{1}x_{5}+5x_{1}x_{6}+7x_{2}x_{3}+3x_{2}x_{4}-5x_{2}x_{5}+3x_{2}x_{6}+7x_{3}^2 \\ 
 & & \quad -2x_{3}x_{4}-x_{3}x_{5}+7x_{3}x_{6}-5x_{4}^2+x_{4}x_{5}-x_{5}^2-x_{5}x_{6} \\ 
q_5 & = & 2x_{1}x_{2}-x_{1}x_{3}+x_{1}x_{5}-7x_{1}x_{6}-x_{2}^2-6x_{2}x_{3}-2x_{2}x_{4}+x_{2}x_{5}-8x_{2}x_{6}-2x_{3}^2-3x_{3}x_{4} \\ 
 & & \quad -5x_{3}x_{5}+2x_{4}^2+3x_{4}x_{6}+x_{5}^2-7x_{5}x_{6}-x_{6}^2 \\ 
q_6 & = & x_{1}^2+3x_{1}x_{2}+5x_{1}x_{3}-3x_{1}x_{4}-2x_{1}x_{5}+2x_{1}x_{6}+2x_{2}^2+4x_{2}x_{3}+x_{2}x_{5}+9x_{2}x_{6} \\ 
 & & \quad -x_{3}^2+2x_{3}x_{4}+4x_{3}x_{5}+9x_{4}^2+6x_{4}x_{5}+5x_{4}x_{6}-x_{5}^2+4x_{6}^2 \\ 
q_7 & = & x_{1}^2+2x_{1}x_{2}+3x_{1}x_{3}-3x_{1}x_{4}+x_{1}x_{6}-6x_{2}^2+3x_{2}x_{3}+9x_{2}x_{4}+x_{2}x_{5}+2x_{2}x_{6}+3x_{3}^2 \\ 
 & & \quad +x_{3}x_{4}+7x_{3}x_{5}-x_{3}x_{6}-6x_{4}^2-2x_{4}x_{5}-6x_{4}x_{6}+3x_{5}^2-x_{5}x_{6}-x_{6}^2 \\ 
q_8 & = & x_{1}^2-2x_{1}x_{2}+x_{1}x_{3}-2x_{1}x_{4}+3x_{1}x_{5}+2x_{1}x_{6}-x_{2}x_{3}-x_{2}x_{4}+3x_{2}x_{5}-8x_{2}x_{6} \\ 
 & & \quad -x_{3}^2+7x_{3}x_{4}+8x_{3}x_{5}+12x_{3}x_{6}+4x_{4}^2-2x_{4}x_{5}+2x_{4}x_{6}+6x_{5}x_{6}-5x_{6}^2 \\ 
q_9 & = & x_{1}x_{2}-6x_{1}x_{3}+6x_{1}x_{4}-2x_{1}x_{5}+6x_{1}x_{6}-x_{2}^2+x_{2}x_{3}-5x_{2}x_{4}-2x_{2}x_{5}-2x_{2}x_{6} \\ 
 & & \quad -7x_{3}^2-3x_{3}x_{4}+5x_{3}x_{5}+x_{3}x_{6}-3x_{4}^2+10x_{4}x_{5}+5x_{4}x_{6}-4x_{5}^2+5x_{5}x_{6}+2x_{6}^2.  
\end{array} $$
\normalsize

The {\tt PointSearch} function (see \S\ref{sec:outline} for references) 
finds a solution
$$ \begin{array}{rcl} 
(x_1 : \ldots :  x_6 ) & = & (7439932626 : -837815413 : -525136075 :  \\
& & \qquad 2262805710 : -3465232629 : -1122238333). \end{array} $$
Mapping back to $C_6$ this point becomes
$$ \begin{pmatrix} x_{11} & x_{12} & x_{13} \\
 x_{21} & x_{22} & x_{23} \end{pmatrix} 
= \begin{pmatrix}
    42664889066 & -110100049465 & -1219600972219 \\
   -157741863480 &  -574453039469 &  114558712088
\end{pmatrix} $$
The $2 \times 2$ minors are the co-ordinates of a point in $C_3(\Q)$
which then maps down to a point $P \in E(\Q)$ of 
canonical height $308.94$.
It is routine to check (using the bounds given in \cite{CPS})
that $P$ is a generator for the Mordell-Weil group modulo torsion.

The final values of $x$ and $y$ for which both 
$x^2 + 809xy + y^2$ and $x^2 -809xy+y^2$ 
are squares may be found on the website \cite{MacLeod}. 
They have $534$ and $537$ decimal digits respectively. 

\subsection{An example of 12-descent}
As described in the introduction, Mark Watkins provided me with
a list of 35 elliptic curves over $\Q$ of analytic rank 2
for which only one generator of the Mordell-Weil group was known.
The Birch--Swinnerton-Dyer conjecture gives an estimate for 
the height of the supposed second generator. The curves
were ordered by this estimated height. For the first 30 curves
on the list, the estimated height was in the range 220 up to 370. The
last 5 were as follows. We list the conductor $N_E$,
the coefficients $a_1, \ldots, a_6$ of a minimal Weierstrass equation,
and the canonical heights of the generators. (The last column was 
computed in hindsight.) 
$$ \begin{array}{rlrr} 
\multicolumn{1}{c}{N_E} & \multicolumn{1}{c}{[a_1,a_2,a_3,a_4,a_6]} 
&  \multicolumn{1}{c}{\widehat{h}(P_1)} 
&  \multicolumn{1}{c}{\widehat{h}(P_2)} \\  \hline
8423178259 & [ 0, -1, 1, -6286122, -6064183289 ] & 17.2636 & 442.070 \\
4817824003 & [ 0, -1, 1, -91969194, -339447383999 ] & 15.4617 & 445.878 \\
4353186907 & [ 1, 1, 0, -14176508, -20550712585 ] & 14.4505 & 488.336 \\
5242805459 & [ 1, 1, 0, -5078887, -4407675042 ] & 2.9643 & 527.301 \\
7800899941 & [ 0, 0, 1, -237882589, -1412186639384 ] & 5.3208 & 642.626 
\end{array} $$

Each curve on the list is the only curve in its isogeny class,
and so in particular has trivial torsion subgroup. 
According to Magma we have $S^{(3)}(E/\Q) \isom (\Z/3\Z)^2$ 
and $S^{(4)}(E/\Q) \isom (\Z/4\Z)^2$. Magma also returns 4~ternary cubics
and 6~pairs of quadrics, representing the inverse pairs of elements
of order $n$ in $S^{(n)}(E/\Q)$ for $n=3,4$. 
Following the method described in \S\ref{sec:outline} we compute
48 different $12$-coverings, each corresponding to an inverse pair
of elements of order $12$ in $S^{(12)}(E/\Q) \isom (\Z/12\Z)^2$.
Since $E(\Q)$ has one generator of small height, we only needed
to search on a few of these before a second independent generator was found.

We give brief details for the last curve on the list. In this
case the first generator is
$$ P_1 = (-2003564/15^2 , -1691/15^3). $$
To compute a second independent generator we start with a 3-covering
$$ C_3 = \left\{ \begin{array}{rcl}
13 x^3 - 87 y^3 - 142 z^3 + 17 x^2 y + 28 x^2 z + 77 x y^2  
\quad \\ 
- \, 187 y^2 z - 35 x z^2 - 30 y z^2 - 118 x y z
 & = & 0 \end{array} \right\} $$
and a 4-covering
$$ C_4 = \left\{ \begin{array}{rcl}
\multicolumn{1}{l}{x_1^2 + 3 x_1 x_2 + 13 x_1 x_4 - 2 x_2^2 - 2 x_2 x_3 \qquad \qquad } \\ 
  + \, x_2 x_4 - 6 x_3^2 - 9 x_3 x_4 + 7 x_4^2 & = & 0 \\
\multicolumn{1}{l}{3 x_1^2 - 6 x_1 x_2 + 5 x_1 x_3 - 14 x_1 x_4 + 8 x_2^2 \qquad \qquad }  \\ 
  - \, 7 x_2 x_3 - 2 x_2 x_4 + 5 x_3^2 + 6 x_4^2 & = & 0
\end{array} \right\}. $$
Following the method described in \S\ref{sec:outline} we compute 54 quadrics
in 12 variables defining a 12-covering $C_{12} \subset \PP(\Mat_{3,4})$.
The coefficients are integers of absolute value at most $40$. After 
minimising at the unique bad prime of $E$, the largest absolute value
was $7$. On this modified curve, {\tt PointSearch} 
(see \S\ref{sec:outline} for references) found a solution
$$ \begin{array}{c} ( -38935814 : 66676907 : 35419393 : 
-17989378 : 14587909 : -9597188 : \\ 
  -41856515 : -6994528 : -103052506 : 12269644 : 
11697462 : 25846956 ) \end{array} $$
Mapping back to $C_{12}$ this point becomes
$$ \begin{pmatrix}
 -585852746652 & -134738830676 &  992806781984 & -476555121265 \\
   -5994121237 &    8026743882 &  -211970353911 &  286395682995 \\
  303306392932 & -167866472332 &  -273061778593 &  215669566507
\end{pmatrix}. $$
The $3 \times 3$ minors are the co-ordinates of a point in $C_4(\Q)$
which then maps down to a point $P =(r/t^2,s/t^3) \in E(\Q)$ of 
canonical height $651.86$ where
$$ \begin{array}{rcl}
t & = & 19114217356093463705777747876066898415631548291608697 \backslash \\ 
  &   & 40922807612824612875940389382477232533975065261036903 \backslash \\ \medskip
  &   & 1136244375962645684728831244647511 \\ 
r & = & 93385419996781156236208893304670769704360761931620474 \backslash \\  
  &   & 91160376652094516256058095438975234936485365750728672 \backslash \\ 
  &   & 93638862617394747880602761519393543195699455909538302 \backslash \\ 
  &   & 59168129312401737073248837456279406678810951156628252 \backslash \\ 
  &   & 40211217008647003170248465787238475381689553329226658 \backslash \\ \medskip
  &   & 862657964535534165 \\ 
s & = & 21189601910515224224247520792578674272370041362778083 \backslash \\  
  &   & 56705954773720391166818153294963600750782215820469113 \backslash \\ 
  &   & 74353930791392149260850703573807892173379962268109766 \backslash \\ 
  &   & 98439592570904852474980215470887488235939468315716611 \backslash \\ 
  &   & 87491555874815362407229178054307290009804071221273367 \backslash \\ 
  &   & 65774805454336495291566121830488793684956520543942634 \backslash \\ 
  &   & 32595140366259647660234205784539280961702449802725098 \backslash \\ 
  &   & 961125300545865563681315860704624955352014647220765212 \\
\end{array} $$ 
A second generator of slightly smaller height is $P_2 = P + P_1$
with $\widehat{h}(P_2) = 642.63.$ (In hindsight we could find $P_2$
directly by starting with different $C_3$ and $C_4$.)
According to Magma the regulator of the subgroup generated by
$P_1$ and $P_2$ is $3415.49$, the non-zero value confirming that 
these points are independent. Again it is routine to check 
(using the bounds given in \cite{CPS}) that $P_1$ and $P_2$
generate the Mordell-Weil group.


\begin{thebibliography}{99}

\frenchspacing
\renewcommand{\baselinestretch}{1}

\bibitem{AKM3P}
S.Y. An, S.Y. Kim, D.C. Marshall, S.H. Marshall,
W.G. McCallum and A.R. Perlis,
Jacobians of genus one curves,
{\em J. Number Theory} 90 (2001), no. 2, 304--315. 

\bibitem{BSD}
B.J. Birch and H.P.F. Swinnerton-Dyer,
Notes on elliptic curves I.  
{\em J. Reine Angew. Math.}  {\bf{212}} (1963), 7--25. 

\bibitem{CaL} 
J.W.S. Cassels,
{\em Lectures on elliptic curves},
LMS Student Texts 24,
Cambridge University Press, Cambridge, 1991. 

\bibitem{CreRed}
J.E. Cremona, 
Reduction of binary cubic and quartic forms,
{\em LMS J. Comput. Math.} {\bf{2}} (1999), 64--94 (electronic).

\bibitem{descsum} 
J.E. Cremona, T.A. Fisher, C. O'Neil, D. Simon and M. Stoll,
{\em Explicit $n$-descent on elliptic curves, I Algebra}, to appear
{\em J. Reine Angew. Math.},
{\em II Geometry}, submitted for publication, 
{\em III Algorithms}, in preparation.

\bibitem{paperIV} 
J.E. Cremona, T.A. Fisher and M. Stoll,
{\em Minimisation and reduction for 3- and 4-coverings of elliptic curves}, 
in preparation.

\bibitem{CPS}
J.E. Cremona, M. Prickett and S. Siksek, 
Height difference bounds for elliptic curves over number fields,
{\em J. Number Theory} {\bf{116}} (2006), no. 1, 42--68.

\bibitem{CS}
J.E. Cremona and M. Stoll, 
Minimal models for 2-coverings of elliptic curves,
{\em LMS J. Comput. Math.} {\bf{5}} (2002), 220--243 (electronic).

\bibitem{Dickson}
L.E. Dickson, 
{\em History of the theory of numbers, Vol. II: Diophantine analysis}, 
Chelsea Publishing Co., New York 1966.

\bibitem{Elk00}
N.D. Elkies, Rational points near curves and small nonzero 
$\vert x\sp 3-y\sp 2\vert$ via lattice reduction, 
{\em Algorithmic number theory} (Leiden, 2000), 33--63, 
Lecture Notes in Comput. Sci., 1838, Springer, Berlin, 2000.

\bibitem{g1inv}
T.A. Fisher, {\em The invariants of a genus one curve}, 
preprint, available at \\ {\tt http://arxiv.org/abs/math/0610318}

\bibitem{Hindry-Silverman}  M. Hindry and J.H. Silverman, 
{\em Diophantine geometry}, Graduate Texts in Mathematics 201, 
Springer-Verlag, New York, 2000. 

\bibitem{Hulek}
K. Hulek, 
{\em Projective geometry of elliptic curves},
Ast\'erisque No. 137 (1986).

\bibitem{MacLeod} 
A.J. MacLeod, Elliptic curves in recreational number theory, website at \\
{\tt http://maths.paisley.ac.uk/allanm/ECRNT/Ecrnt.htm}

\bibitem{Magma} {\sf MAGMA} is described in
W. Bosma, J. Cannon and C. Playoust,
The Magma algebra system I: The user language, 
{\em J. Symbolic Comput.} 24, 235--265 (1997).
The Magma home page is at {\tt http://magma.maths.usyd.edu.au/magma/}

\bibitem{Mumford}
D. Mumford, {\em Abelian varieties}, Oxford University Press, 1970.

\bibitem{Cathy} 
C. O'Neil, The period-index obstruction for elliptic curves,
{\em  J. Number Theory}  95  (2002),  no. 2, 329--339. 

\bibitem{Siksek}
S. Siksek, {\em Descent on curve of genus $1$}, 
PhD thesis, University of Exeter, 1995. 

\bibitem{Silverman}
J.H. Silverman,
{\em The arithmetic of elliptic curves},
Graduate Texts in Mathematics 106, Springer-Verlag, New York, 1992. 
                                                         
\bibitem{Stamminger}
S. Stamminger,
{\em Explicit 8-descent on elliptic curves},
PhD thesis, International University Bremen, 2005.

\bibitem{SteinWatkins}
W.A. Stein and M. Watkins, 
A database of elliptic curves---first report, 
{\em Algorithmic number theory} (Sydney, 2002), 267--275,
Lecture Notes in Comput. Sci., 2369, Springer, Berlin, 2002.

\bibitem{Watkins}
M. Watkins, {\em Searching for points $p$-adically}, notes available
from \\ {\tt http://www.maths.bris.ac.uk/$\sim$mamjw/papers/}

\bibitem{Weil}
A. Weil, 
Remarques sur un m\'emoire d'Hermite,
{\em Arch. Math.} 5, (1954). 197--202.

\bibitem{Womack}
T. Womack,
{\em Explicit descent on elliptic curves},
PhD thesis, University of Nottingham, 2003.

\bibitem{Zarhin}
Ju. G. Zarhin, 
Noncommutative cohomology and Mumford groups,
{\em Math. Notes} {\bf{15}} (1974), 241--244.

\end{thebibliography}
\end{document}